%%%%%%%%%%%%%%%%%%%%%%%%%%%%%%%%%%%%%%%%%%%%%%%%%%%%%%%%%%%%%%%%%%%%%%%%%
%									%
%	Duality for Knizhnik-Zamolodchikov and Dynamical Equations	%
%									%
%	V.Tarasov  and  A.Varchenko					%
%									%
%	Preprint (2001), 12 pages, amstex 2.2				%
%									%
%	Acta Applicandae Mathematicae 73 (2002)	no. 1-2, 141--154	%
%	Proceedings of The 2000 Twente conference on Lie groups		%
%									%
%%%%%%%%%%%%%%%%%%%%%%%%%%%%%%%%%%%%%%%%%%%%%%%%%%%%%%%%%%%%%%%%%%%%%%%%%

\mag 1200

\input amstex

\expandafter\ifx\csname dual.def\endcsname\relax \else\endinput\fi
\expandafter\edef\csname dual.def\endcsname{%
 \catcode`\noexpand\@=\the\catcode`\@\space}

\let\atbefore @

\catcode`\@=11

\overfullrule\z@

\def\PaperA4{\hsize 6.25truein \vsize 9.63truein}

 %% American paper 8.5x11

\def\foliorm{\ifMag\eightrm\else\ninerm\fi}

\let\@ft@\expandafter \let\@tb@f@\atbefore

\newif\ifMag
\def\Magset{\ifnum\mag>\@m\Magtrue\fi}
\Magset

\newif\ifUS

\newdimen\p@@ \p@@\p@
\def\m@ths@r{\ifnum\mathsurround=\z@\z@\else\maths@r\fi}
\def\maths@r{1.6\p@@} \def\mathsurzero{\def\maths@r{\z@}}

\mathsurround\maths@r
\font\Brm=cmr12 \font\Bbf=cmbx12 \font\Bit=cmti12 \font\ssf=cmss10
\font\Bsl=cmsl10 scaled 1200 \font\Bmmi=cmmi10 scaled 1200
\font\BBf=cmbx12 scaled 1200 \font\BMmi=cmmi10 scaled 1440

\def\atletter{\edef\atrestore{\catcode`\noexpand\@=\the\catcode`\@\space}
 \catcode`\@=11}

\newread\@ux \newwrite\@@x \newwrite\@@cd
\let\@np@@\input
\def\@np@t#1{\openin\@ux#1\relax\ifeof\@ux\else\closein\@ux\relax\@np@@ #1\fi}
\def\input#1 {\openin\@ux#1\relax\ifeof\@ux\wrs@x{No file #1}\else
 \closein\@ux\relax\@np@@ #1\fi}
\def\Input#1 {\relax} %% Do not remove the space after #1

\def\wr@@x#1{} \def\wrs@x{\immediate\write\sixt@@n}

\def\readldf{\@np@t{\jobname.ldf}}
\def\writeldf{\def\wr@@x{\immediate\write\@@x}\def\wr@x@{\write\@@x}
 \def\cl@selbl{\wr@@x{\string\def\string\nextpage{\the\pageno}}%
 \wr@@x{\string\endinput}\immediate\closeout\@@x}
 \immediate\openout\@@x\jobname.ldf}
\let\cl@selbl\relax

\def\nextpage{1}

\def\tod@y{\ifcase\month\or
 January\or February\or March\or April\or May\or June\or July\or
 August\or September\or October\or November\or December\fi\space\,
\number\day,\space\,\number\year}

\newcount\c@time
\def\h@@r{hh}\def\m@n@te{mm}
\def\wh@tt@me{\c@time\time\divide\c@time 60\edef\h@@r{\number\c@time}%
 \multiply\c@time -60\advance\c@time\time\edef
 \m@n@te{\ifnum\c@time<10 0\fi\number\c@time}}
\def\t@me{\h@@r\/{\rm:}\m@n@te} \let\whattime\wh@tt@me
\def\today{\tod@y\wr@@x{\string\todaydef{\tod@y}}}
\def\nowtime{\t@me{\let\/\ic@\wr@@x{\string\nowtimedef{\t@me}}}}
\def\todaydef#1{} \def\nowtimedef#1{}

\def\em#1{{\it #1\/}} \def\emph#1{{\sl #1\/}}

\def\fitem#1{\par\setbox\z@\hbox{#1}\hangindent\wd\z@
 \hglue-2\parindent\kern\wd\z@\indent\llap{#1}\ignore}

\def\itemflat#1{\par\setbox\z@\hbox{\rm #1\enspace}\hang\ifnum\wd\z@>\parindent
 \noindent\unhbox\z@\ignore\else\textindent{\rm#1}\fi}

\newcount\itemlet
\def\newbi{\itemlet 96} \newbi
\def\bitem{\gad\itemlet \par\hangindent1.5\parindent
 \hglue-.5\parindent\textindent{\rm\rlap{\char\the\itemlet}\hp{b})}}

\newcount\itemrm

\def\iitem{\gad\itemrm \par\hangindent1.5\parindent
 \hglue-.5\parindent\textindent{\rm\hp{v}\llap{\romannumeral\the\itemrm})}}

\newcount\itemar

\def\iitema{\gad\itemrm \par\hangindent1.5\parindent
 \hglue-.5\parindent\textindent{\rm\hp{0}\llap{\the\itemrm}.}}

\def\center{\par\begingroup\leftskip\z@ plus \hsize \rightskip\leftskip
 \parindent\z@\parfillskip\z@skip \def\\{\unskip\break}}
\def\endcenter{\endgraf\endgroup}

\def\Abstract{\begingroup\narrower\nt{\bf Abstract.}\enspace\ignore}
\def\endAbs{\endgraf\endgroup}

\let\b@gr@@\begingroup \let\B@gr@@\begingroup
\def\b@gr@{\b@gr@@\let\b@gr@@\undefined}
\def\B@gr@{\B@gr@@\let\B@gr@@\undefined}

\def\@fn@xt#1#2#3{\let\@ch@r=#1\def\n@xt{\ifx\t@st@\@ch@r
 \def\n@@xt{#2}\else\def\n@@xt{#3}\fi\n@@xt}\futurelet\t@st@\n@xt}

\def\@fwd@@#1#2#3{\setbox\z@\hbox{#1}\ifdim\wd\z@>\z@#2\else#3\fi}
\def\s@twd@#1#2{\setbox\z@\hbox{#2}#1\wd\z@}

\def\r@st@re#1{\let#1\s@v@} \def\s@v@d@f{\let\s@v@}

\def\p@sk@p#1#2{\par\skip@#2\relax\ifdim\lastskip<\skip@\relax\removelastskip
 \ifnum#1=\z@\else\penalty#1\relax\fi\vskip\skip@
 \else\ifnum#1=\z@\else\penalty#1\relax\fi\fi}
\def\sk@@p#1{\par\skip@#1\relax\ifdim\lastskip<\skip@\relax\removelastskip
 \vskip\skip@\fi}

\newbox\p@b@ld
\def\poorbold#1{\setbox\p@b@ld\hbox{#1}\kern-.01em\copy\p@b@ld\kern-\wd\p@b@ld
 \kern.02em\copy\p@b@ld\kern-\wd\p@b@ld\kern-.012em\raise.02em\box\p@b@ld}

\ifx\plainfootnote\undefined \let\plainfootnote\footnote \fi

\let\s@v@\proclaim \let\proclaim\relax
\def\r@R@fs#1{\let#1\s@R@fs} \let\s@R@fs\Refs \let\Refs\relax
\def\r@endd@#1{\let#1\s@endd@} \let\s@endd@\enddocument
\let\bye\relax

\def\myR@fs{\@fn@xt[\m@R@f@\m@R@fs} \def\m@R@fs{\@fn@xt*\m@r@f@@\m@R@f@@}
\def\m@R@f@@{\m@R@f@[References]} \def\m@r@f@@*{\m@R@f@[]}

\def\Twelvepoint{\twelvepoint \let\Bbf\BBf \let\Bmmi\BMmi
\font\Brm=cmr12 scaled 1200 \font\Bit=cmti12 scaled 1200
\font\ssf=cmss10 scaled 1200 \font\Bsl=cmsl10 scaled 1440
\font\BBf=cmbx12 scaled 1440 \font\BMmi=cmmi10 scaled 1728}

\newdimen\b@gsize

\newdimen\r@f@nd \newbox\r@f@b@x \newbox\adjb@x
\newbox\p@nct@ \newbox\k@yb@x \newcount\rcount
\newbox\b@b@x \newbox\p@p@rb@x \newbox\j@@rb@x \newbox\y@@rb@x
\newbox\v@lb@x \newbox\is@b@x \newbox\p@g@b@x \newif\ifp@g@ \newif\ifp@g@s
\newbox\inb@@kb@x \newbox\b@@kb@x \newbox\p@blb@x \newbox\p@bl@db@x
\newbox\ed@b@x \newif\ifed@ \newif\ifed@s \newif\if@fl@b \newif\if@fn@m
\newbox\p@p@nf@b@x \newbox\inf@b@x \newbox\b@@nf@b@x
\newtoks\@dd@p@n \newtoks\@ddt@ks

\newif\ifp@gen@

\def\p@@nt{.\kern.3em} \let\point\p@@nt

\let\proheadfont\bf \let\probodyfont\sl \let\demofont\it

\headline={\hfil}
\footline={\ifp@gen@\ifnum\pageno=\z@\else\hfil\foliorm\folio\fi\else
 \ifnum\pageno=\z@\hfil\foliorm\folio\fi\fi\hfil\global\p@gen@true}
\parindent1pc

\font@\tensmc=cmcsc10
\font@\sevenex=cmex7
\font@\sevenit=cmti7
\font@\eightrm=cmr8
\font@\sixrm=cmr6
\font@\eighti=cmmi8 \skewchar\eighti='177
\font@\sixi=cmmi6 \skewchar\sixi='177
\font@\eightsy=cmsy8 \skewchar\eightsy='60
\font@\sixsy=cmsy6 \skewchar\sixsy='60
\font@\eightex=cmex8
\font@\eightbf=cmbx8
\font@\sixbf=cmbx6
\font@\eightit=cmti8
\font@\eightsl=cmsl8
\font@\eightsmc=cmcsc8
\font@\eighttt=cmtt8
\font@\ninerm=cmr9
\font@\ninei=cmmi9 \skewchar\ninei='177
\font@\ninesy=cmsy9 \skewchar\ninesy='60
\font@\nineex=cmex9
\font@\ninebf=cmbx9
\font@\nineit=cmti9
\font@\ninesl=cmsl9
\font@\ninesmc=cmcsc9
\font@\ninemsa=msam9
\font@\ninemsb=msbm9
\font@\nineeufm=eufm9
\font@\eightmsa=msam8
\font@\eightmsb=msbm8
\font@\eighteufm=eufm8
\font@\sixmsa=msam6
\font@\sixmsb=msbm6
\font@\sixeufm=eufm6

\loadmsam\loadmsbm\loadeufm
\input amssym.tex

\def\footnoterule{\kern-3\p@\hrule width5pc\kern 2.6\p@}
\def\m@k@foot#1{\insert\footins
 {\interlinepenalty\interfootnotelinepenalty
 \ifMag\eightpoint\else\ninepoint\fi
 \splittopskip\ht\strutbox\splitmaxdepth\dp\strutbox
 \floatingpenalty\@MM\leftskip\z@\rightskip\z@
 \spaceskip\z@\xspaceskip\z@
 \leavevmode\footstrut\ignore#1\unskip\lower\dp\strutbox
 \vbox to\dp\strutbox{}}}
\def\ftext#1{\m@k@foot{\vsk-.8>\nt #1}}
\def\pr@cl@@m#1{\p@sk@p{-100}\medskipamount
 \def\endproclaim{\endgroup\p@sk@p{55}\medskipamount}\begingroup
 \nt\ignore\proheadfont#1\unskip.\enspace\probodyfont\ignore}
\outer\def\proclaim{\pr@cl@@m} \s@v@d@f\proclaim \let\proclaim\relax
\def\demo#1{\sk@@p\medskipamount\nt{\ignore\demofont#1\unskip.}\enspace
 \ignore}
\def\enddemo{\sk@@p\medskipamount}

\def\cite#1{{\rm[#1]}} 
 \def\Refs#1#2{\relax}

\def\big@#1#2{{\hbox{$\left#2\vcenter to#1\b@gsize{}%
 \right.\nulldelimiterspace\z@\m@th$}}}
\def\big{\big@\@ne}
\def\Big{\big@{1.5}}
\def\bigg{\big@\tw@}
\def\Bigg{\big@{2.5}}
\normallineskiplimit\p@

\def\tenpoint{\p@@\p@ \normallineskiplimit\p@@
 \mathsurround\m@ths@r \normalbaselineskip12\p@@
 \abovedisplayskip12\p@@ plus3\p@@ minus9\p@@
 \belowdisplayskip\abovedisplayskip
 \abovedisplayshortskip\z@ plus3\p@@
 \belowdisplayshortskip7\p@@ plus3\p@@ minus4\p@@
 \textonlyfont@\rm\tenrm \textonlyfont@\it\tenit
 \textonlyfont@\sl\tensl \textonlyfont@\bf\tenbf
 \textonlyfont@\smc\tensmc \textonlyfont@\tt\tentt
 \ifsyntax@ \def\big##1{{\hbox{$\left##1\right.$}}}%
  \let\Big\big \let\bigg\big \let\Bigg\big
 \else
  \textfont\z@\tenrm \scriptfont\z@\sevenrm \scriptscriptfont\z@\fiverm
  \textfont\@ne\teni \scriptfont\@ne\seveni \scriptscriptfont\@ne\fivei
  \textfont\tw@\tensy \scriptfont\tw@\sevensy \scriptscriptfont\tw@\fivesy
  \textfont\thr@@\tenex \scriptfont\thr@@\sevenex
	\scriptscriptfont\thr@@\sevenex
  \textfont\itfam\tenit \scriptfont\itfam\sevenit
	\scriptscriptfont\itfam\sevenit
  \textfont\bffam\tenbf \scriptfont\bffam\sevenbf
	\scriptscriptfont\bffam\fivebf
  \textfont\msafam\tenmsa \scriptfont\msafam\sevenmsa
	\scriptscriptfont\msafam\fivemsa
  \textfont\msbfam\tenmsb \scriptfont\msbfam\sevenmsb
	\scriptscriptfont\msbfam\fivemsb
  \textfont\eufmfam\teneufm \scriptfont\eufmfam\seveneufm
	\scriptscriptfont\eufmfam\fiveeufm
  \setbox\strutbox\hbox{\vrule height8.5\p@@ depth3.5\p@@ width\z@}%
  \setbox\strutbox@\hbox{\lower.5\normallineskiplimit\vbox{%
	\kern-\normallineskiplimit\copy\strutbox}}%
   \setbox\z@\vbox{\hbox{$($}\kern\z@}\b@gsize1.2\ht\z@
  \fi
  \normalbaselines\rm\dotsspace@1.5mu\ex@.2326ex\jot3\ex@}

\def\eightpoint{\p@@.8\p@ \normallineskiplimit\p@@
 \mathsurround\m@ths@r \normalbaselineskip10\p@
 \abovedisplayskip10\p@ plus2.4\p@ minus7.2\p@
 \belowdisplayskip\abovedisplayskip
 \abovedisplayshortskip\z@ plus3\p@@
 \belowdisplayshortskip7\p@@ plus3\p@@ minus4\p@@
 \textonlyfont@\rm\eightrm \textonlyfont@\it\eightit
 \textonlyfont@\sl\eightsl \textonlyfont@\bf\eightbf
 \textonlyfont@\smc\eightsmc \textonlyfont@\tt\eighttt
 \ifsyntax@\def\big##1{{\hbox{$\left##1\right.$}}}%
  \let\Big\big \let\bigg\big \let\Bigg\big
 \else
  \textfont\z@\eightrm \scriptfont\z@\sixrm \scriptscriptfont\z@\fiverm
  \textfont\@ne\eighti \scriptfont\@ne\sixi \scriptscriptfont\@ne\fivei
  \textfont\tw@\eightsy \scriptfont\tw@\sixsy \scriptscriptfont\tw@\fivesy
  \textfont\thr@@\eightex \scriptfont\thr@@\sevenex
	\scriptscriptfont\thr@@\sevenex
  \textfont\itfam\eightit \scriptfont\itfam\sevenit
	\scriptscriptfont\itfam\sevenit
  \textfont\bffam\eightbf \scriptfont\bffam\sixbf
	\scriptscriptfont\bffam\fivebf
  \textfont\msafam\eightmsa \scriptfont\msafam\sixmsa
	\scriptscriptfont\msafam\fivemsa
  \textfont\msbfam\eightmsb \scriptfont\msbfam\sixmsb
	\scriptscriptfont\msbfam\fivemsb
  \textfont\eufmfam\eighteufm \scriptfont\eufmfam\sixeufm
	\scriptscriptfont\eufmfam\fiveeufm
 \setbox\strutbox\hbox{\vrule height7\p@ depth3\p@ width\z@}%
 \setbox\strutbox@\hbox{\raise.5\normallineskiplimit\vbox{%
   \kern-\normallineskiplimit\copy\strutbox}}%
 \setbox\z@\vbox{\hbox{$($}\kern\z@}\b@gsize1.2\ht\z@
 \fi
 \normalbaselines\eightrm\dotsspace@1.5mu\ex@.2326ex\jot3\ex@}

\def\ninepoint{\p@@.9\p@ \normallineskiplimit\p@@
 \mathsurround\m@ths@r \normalbaselineskip11\p@
 \abovedisplayskip11\p@ plus2.7\p@ minus8.1\p@
 \belowdisplayskip\abovedisplayskip
 \abovedisplayshortskip\z@ plus3\p@@
 \belowdisplayshortskip7\p@@ plus3\p@@ minus4\p@@
 \textonlyfont@\rm\ninerm \textonlyfont@\it\nineit
 \textonlyfont@\sl\ninesl \textonlyfont@\bf\ninebf
 \textonlyfont@\smc\ninesmc \textonlyfont@\tt\ninett
 \ifsyntax@ \def\big##1{{\hbox{$\left##1\right.$}}}%
  \let\Big\big \let\bigg\big \let\Bigg\big
 \else
  \textfont\z@\ninerm \scriptfont\z@\sevenrm \scriptscriptfont\z@\fiverm
  \textfont\@ne\ninei \scriptfont\@ne\seveni \scriptscriptfont\@ne\fivei
  \textfont\tw@\ninesy \scriptfont\tw@\sevensy \scriptscriptfont\tw@\fivesy
  \textfont\thr@@\nineex \scriptfont\thr@@\sevenex
	\scriptscriptfont\thr@@\sevenex
  \textfont\itfam\nineit \scriptfont\itfam\sevenit
	\scriptscriptfont\itfam\sevenit
  \textfont\bffam\ninebf \scriptfont\bffam\sevenbf
	\scriptscriptfont\bffam\fivebf
  \textfont\msafam\ninemsa \scriptfont\msafam\sevenmsa
	\scriptscriptfont\msafam\fivemsa
  \textfont\msbfam\ninemsb \scriptfont\msbfam\sevenmsb
	\scriptscriptfont\msbfam\fivemsb
  \textfont\eufmfam\nineeufm \scriptfont\eufmfam\seveneufm
	\scriptscriptfont\eufmfam\fiveeufm
  \setbox\strutbox\hbox{\vrule height8.5\p@@ depth3.5\p@@ width\z@}%
  \setbox\strutbox@\hbox{\lower.5\normallineskiplimit\vbox{%
	\kern-\normallineskiplimit\copy\strutbox}}%
   \setbox\z@\vbox{\hbox{$($}\kern\z@}\b@gsize1.2\ht\z@
  \fi
  \normalbaselines\rm\dotsspace@1.5mu\ex@.2326ex\jot3\ex@}

\font@\twelverm=cmr10 scaled 1200
\font@\twelveit=cmti10 scaled 1200
\font@\twelvesl=cmsl10 scaled 1200
\font@\twelvebf=cmbx10 scaled 1200
\font@\twelvesmc=cmcsc10 scaled 1200
\font@\twelvett=cmtt10 scaled 1200
\font@\twelvei=cmmi10 scaled 1200 \skewchar\twelvei='177
\font@\twelvesy=cmsy10 scaled 1200 \skewchar\twelvesy='60
\font@\twelveex=cmex10 scaled 1200
\font@\twelvemsa=msam10 scaled 1200
\font@\twelvemsb=msbm10 scaled 1200
\font@\twelveeufm=eufm10 scaled 1200

\def\twelvepoint{\p@@1.2\p@ \normallineskiplimit\p@@
 \mathsurround\m@ths@r \normalbaselineskip12\p@@
 \abovedisplayskip12\p@@ plus3\p@@ minus9\p@@
 \belowdisplayskip\abovedisplayskip
 \abovedisplayshortskip\z@ plus3\p@@
 \belowdisplayshortskip7\p@@ plus3\p@@ minus4\p@@
 \textonlyfont@\rm\twelverm \textonlyfont@\it\twelveit
 \textonlyfont@\sl\twelvesl \textonlyfont@\bf\twelvebf
 \textonlyfont@\smc\twelvesmc \textonlyfont@\tt\twelvett
 \ifsyntax@ \def\big##1{{\hbox{$\left##1\right.$}}}%
  \let\Big\big \let\bigg\big \let\Bigg\big
 \else
  \textfont\z@\twelverm \scriptfont\z@\eightrm \scriptscriptfont\z@\sixrm
  \textfont\@ne\twelvei \scriptfont\@ne\eighti \scriptscriptfont\@ne\sixi
  \textfont\tw@\twelvesy \scriptfont\tw@\eightsy \scriptscriptfont\tw@\sixsy
  \textfont\thr@@\twelveex \scriptfont\thr@@\eightex
	\scriptscriptfont\thr@@\sevenex
  \textfont\itfam\twelveit \scriptfont\itfam\eightit
	\scriptscriptfont\itfam\sevenit
  \textfont\bffam\twelvebf \scriptfont\bffam\eightbf
	\scriptscriptfont\bffam\sixbf
  \textfont\msafam\twelvemsa \scriptfont\msafam\eightmsa
	\scriptscriptfont\msafam\sixmsa
  \textfont\msbfam\twelvemsb \scriptfont\msbfam\eightmsb
	\scriptscriptfont\msbfam\sixmsb
  \textfont\eufmfam\twelveeufm \scriptfont\eufmfam\eighteufm
	\scriptscriptfont\eufmfam\sixeufm
  \setbox\strutbox\hbox{\vrule height8.5\p@@ depth3.5\p@@ width\z@}%
  \setbox\strutbox@\hbox{\lower.5\normallineskiplimit\vbox{%
	\kern-\normallineskiplimit\copy\strutbox}}%
  \setbox\z@\vbox{\hbox{$($}\kern\z@}\b@gsize1.2\ht\z@
  \fi
  \normalbaselines\rm\dotsspace@1.5mu\ex@.2326ex\jot3\ex@}

\font@\twelvetrm=cmr10 at 12truept
\font@\twelvetit=cmti10 at 12truept
\font@\twelvetsl=cmsl10 at 12truept
\font@\twelvetbf=cmbx10 at 12truept
\font@\twelvetsmc=cmcsc10 at 12truept
\font@\twelvettt=cmtt10 at 12truept
\font@\twelveti=cmmi10 at 12truept \skewchar\twelveti='177
\font@\twelvetsy=cmsy10 at 12truept \skewchar\twelvetsy='60
\font@\twelvetex=cmex10 at 12truept
\font@\twelvetmsa=msam10 at 12truept
\font@\twelvetmsb=msbm10 at 12truept
\font@\twelveteufm=eufm10 at 12truept

\def\twelvetruepoint{\p@@1.2truept \normallineskiplimit\p@@
 \mathsurround\m@ths@r \normalbaselineskip12\p@@
 \abovedisplayskip12\p@@ plus3\p@@ minus9\p@@
 \belowdisplayskip\abovedisplayskip
 \abovedisplayshortskip\z@ plus3\p@@
 \belowdisplayshortskip7\p@@ plus3\p@@ minus4\p@@
 \textonlyfont@\rm\twelvetrm \textonlyfont@\it\twelvetit
 \textonlyfont@\sl\twelvetsl \textonlyfont@\bf\twelvetbf
 \textonlyfont@\smc\twelvetsmc \textonlyfont@\tt\twelvettt
 \ifsyntax@ \def\big##1{{\hbox{$\left##1\right.$}}}%
  \let\Big\big \let\bigg\big \let\Bigg\big
 \else
  \textfont\z@\twelvetrm \scriptfont\z@\eightrm \scriptscriptfont\z@\sixrm
  \textfont\@ne\twelveti \scriptfont\@ne\eighti \scriptscriptfont\@ne\sixi
  \textfont\tw@\twelvetsy \scriptfont\tw@\eightsy \scriptscriptfont\tw@\sixsy
  \textfont\thr@@\twelvetex \scriptfont\thr@@\eightex
	\scriptscriptfont\thr@@\sevenex
  \textfont\itfam\twelvetit \scriptfont\itfam\eightit
	\scriptscriptfont\itfam\sevenit
  \textfont\bffam\twelvetbf \scriptfont\bffam\eightbf
	\scriptscriptfont\bffam\sixbf
  \textfont\msafam\twelvetmsa \scriptfont\msafam\eightmsa
	\scriptscriptfont\msafam\sixmsa
  \textfont\msbfam\twelvetmsb \scriptfont\msbfam\eightmsb
	\scriptscriptfont\msbfam\sixmsb
  \textfont\eufmfam\twelveteufm \scriptfont\eufmfam\eighteufm
	\scriptscriptfont\eufmfam\sixeufm
  \setbox\strutbox\hbox{\vrule height8.5\p@@ depth3.5\p@@ width\z@}%
  \setbox\strutbox@\hbox{\lower.5\normallineskiplimit\vbox{%
	\kern-\normallineskiplimit\copy\strutbox}}%
  \setbox\z@\vbox{\hbox{$($}\kern\z@}\b@gsize1.2\ht\z@
  \fi
  \normalbaselines\rm\dotsspace@1.5mu\ex@.2326ex\jot3\ex@}

\font@\elevenrm=cmr10 scaled 1095
\font@\elevenit=cmti10 scaled 1095
\font@\elevensl=cmsl10 scaled 1095
\font@\elevenbf=cmbx10 scaled 1095
\font@\elevensmc=cmcsc10 scaled 1095
\font@\eleventt=cmtt10 scaled 1095
\font@\eleveni=cmmi10 scaled 1095 \skewchar\eleveni='177
\font@\elevensy=cmsy10 scaled 1095 \skewchar\elevensy='60
\font@\elevenex=cmex10 scaled 1095
\font@\elevenmsa=msam10 scaled 1095
\font@\elevenmsb=msbm10 scaled 1095
\font@\eleveneufm=eufm10 scaled 1095

\def\elevenpoint{\p@@1.1\p@ \normallineskiplimit\p@@
 \mathsurround\m@ths@r \normalbaselineskip12\p@@
 \abovedisplayskip12\p@@ plus3\p@@ minus9\p@@
 \belowdisplayskip\abovedisplayskip
 \abovedisplayshortskip\z@ plus3\p@@
 \belowdisplayshortskip7\p@@ plus3\p@@ minus4\p@@
 \textonlyfont@\rm\elevenrm \textonlyfont@\it\elevenit
 \textonlyfont@\sl\elevensl \textonlyfont@\bf\elevenbf
 \textonlyfont@\smc\elevensmc \textonlyfont@\tt\eleventt
 \ifsyntax@ \def\big##1{{\hbox{$\left##1\right.$}}}%
  \let\Big\big \let\bigg\big \let\Bigg\big
 \else
  \textfont\z@\elevenrm \scriptfont\z@\eightrm \scriptscriptfont\z@\sixrm
  \textfont\@ne\eleveni \scriptfont\@ne\eighti \scriptscriptfont\@ne\sixi
  \textfont\tw@\elevensy \scriptfont\tw@\eightsy \scriptscriptfont\tw@\sixsy
  \textfont\thr@@\elevenex \scriptfont\thr@@\eightex
	\scriptscriptfont\thr@@\sevenex
  \textfont\itfam\elevenit \scriptfont\itfam\eightit
	\scriptscriptfont\itfam\sevenit
  \textfont\bffam\elevenbf \scriptfont\bffam\eightbf
	\scriptscriptfont\bffam\sixbf
  \textfont\msafam\elevenmsa \scriptfont\msafam\eightmsa
	\scriptscriptfont\msafam\sixmsa
  \textfont\msbfam\elevenmsb \scriptfont\msbfam\eightmsb
	\scriptscriptfont\msbfam\sixmsb
  \textfont\eufmfam\eleveneufm \scriptfont\eufmfam\eighteufm
	\scriptscriptfont\eufmfam\sixeufm
  \setbox\strutbox\hbox{\vrule height8.5\p@@ depth3.5\p@@ width\z@}%
  \setbox\strutbox@\hbox{\lower.5\normallineskiplimit\vbox{%
	\kern-\normallineskiplimit\copy\strutbox}}%
  \setbox\z@\vbox{\hbox{$($}\kern\z@}\b@gsize1.2\ht\z@
  \fi
  \normalbaselines\rm\dotsspace@1.5mu\ex@.2326ex\jot3\ex@}

\def\m@R@f@[#1]{\mathsurzero{%\let\{\relax
 \s@ct{}{#1}}\wr@@c{\string\Refcd{#1}{\the\pageno}}\B@gr@
 \frenchspacing\rcount\z@\refkey{\k@yf@nt[##1]}\refno{\k@yf@nt[##1]}%
 \widest{AZ}\keyright\let\Key\key\let\refin\relax}
\def\widest#1{\s@twd@\r@f@nd{\r@fk@y{\k@yf@nt#1}\enspace}}
\def\widestno#1{\s@twd@\r@f@nd{\r@fn@{\k@yf@nt#1}\enspace}}
\def\widestlabel#1{\s@twd@\r@f@nd{\k@yf@nt#1\enspace}}
\def\refkey{\def\r@fk@y##1} \def\refno{\def\r@fn@##1}
\def\keyright{\def\r@fit@m{\hang\textindent}}
\def\keyflat{\def\r@fit@m##1{\setbox\z@\hbox{##1\enspace}\hang\noindent
 \ifnum\wd\z@<\parindent\indent\hglue-\wd\z@\fi\unhbox\z@}}

\def\R@fb@x{\global\setbox\r@f@b@x} \def\K@yb@x{\global\setbox\k@yb@x}
\def\ref{\par\b@gr@\r@ff@nt\R@fb@x\box\voidb@x\K@yb@x\box\voidb@x
 \@fn@mfalse\@fl@bfalse\b@g@nr@f}
\def\c@nc@t#1{\setbox\z@\lastbox
 \setbox\adjb@x\hbox{\unhbox\adjb@x\unhbox\z@\unskip\unskip\unpenalty#1}}
\def\adjust#1{\relax\ifmmode\penalty-\@M\null\hfil$\clubpenalty\z@
 \widowpenalty\z@\interlinepenalty\z@\offinterlineskip\endgraf
 \setbox\z@\lastbox\unskip\unpenalty\c@nc@t{#1}\nt$\hfil\penalty-\@M
 \else\endgraf\c@nc@t{#1}\nt\fi}
\def\adjustnext#1{\P@nct\hbox{#1}\ignore}
\def\adjustend#1{\def\@djp@{#1}\ignore}
\def\addtoks#1{\global\@ddt@ks{#1}\ignore}
\def\addnext#1{\global\@dd@p@n{#1}\ignore}

\def\cl@s@{\adjust{\@djp@}\endgraf\setbox\z@\lastbox
 \global\setbox\@ne\hbox{\unhbox\adjb@x\ifvoid\z@\else\unhbox\z@\unskip\unskip
 \unpenalty\fi}\egroup\ifnum\c@rr@nt=\k@yb@x\global\fi
 \setbox\c@rr@nt\hbox{\unhbox\@ne\box\p@nct@}\P@nct\null
 \the\@ddt@ks\global\@ddt@ks{}}
\def\@p@n#1{\def\c@rr@nt{#1}\setbox\c@rr@nt\vbox\bgroup\let\@djp@\relax
 \hsize\maxdimen\nt\the\@dd@p@n\global\@dd@p@n{}}
\def\b@g@nr@f{\bgroup\@p@n\z@}
\def\key{\cl@s@\ifvoid\k@yb@x\@p@n\k@yb@x\k@yf@nt\else\@p@n\z@\fi}
\def\label{\cl@s@\ifvoid\k@yb@x\global\@fl@btrue\@p@n\k@yb@x\k@yf@nt\else
 \@p@n\z@\fi}
\def\no{\cl@s@\ifvoid\k@yb@x\gad\rcount\global\@fn@mtrue
 \K@yb@x\hbox{\k@yf@nt\the\rcount}\fi\@p@n\z@}
\def\labelno{\cl@s@\ifvoid\k@yb@x\gad\rcount\@fl@btrue
 \@p@n\k@yb@x\k@yf@nt\the\rcount\else\@p@n\z@\fi}
\def\by{\cl@s@\@p@n\b@b@x} \def\paper{\cl@s@\@p@n\p@p@rb@x\p@p@rf@nt\ignore}
\def\jour{\cl@s@\@p@n\j@@rb@x} \def\yr{\cl@s@\@p@n\y@@rb@x}
\def\vol{\cl@s@\@p@n\v@lb@x\v@lf@nt\ignore}
\def\issue{\cl@s@\@p@n\is@b@x\iss@f@nt\ignore}
\def\page{\cl@s@\ifp@g@s\@p@n\z@\else\p@g@true\@p@n\p@g@b@x\fi}
\def\pages{\cl@s@\ifp@g@\@p@n\z@\else\p@g@strue\@p@n\p@g@b@x\fi}
\def\inbook{\cl@s@\@p@n\inb@@kb@x}
\def\book{\cl@s@\@p@n\b@@kb@x\b@@kf@nt\ignore}
\def\publ{\cl@s@\@p@n\p@blb@x} \def\publaddr{\cl@s@\@p@n\p@bl@db@x}
\def\ed{\cl@s@\ifed@s\@p@n\z@\else\ed@true\@p@n\ed@b@x\fi}
\def\eds{\cl@s@\ifed@\@p@n\z@\else\ed@strue\@p@n\ed@b@x\fi}
\def\info{\cl@s@\@p@n\inf@b@x} \def\paperinfo{\cl@s@\@p@n\p@p@nf@b@x}
\def\bookinfo{\cl@s@\@p@n\b@@nf@b@x} 
\def\P@nct{\global\setbox\p@nct@} \def\nopunct{\P@nct\box\voidb@x}
\def\p@@@t#1#2{\ifvoid\p@nct@\else#1\unhbox\p@nct@#2\fi}
\def\sp@@{\penalty-50 \space\hskip\z@ plus.1em}
\def\c@mm@{\p@@@t,\sp@@} \def\sp@c@{\p@@@t\empty\sp@@}
\def\p@tb@x#1#2{\ifvoid#1\else#2\@nb@x#1\fi}
\def\@nb@x#1{\unhbox#1\P@nct\lastbox}
\def\endr@f@{\cl@s@\nopunct
 \R@fb@x\hbox{\unhbox\r@f@b@x \p@tb@x\b@b@x\empty
 \ifvoid\j@@rb@x\ifvoid\inb@@kb@x\ifvoid\p@p@rb@x\ifvoid\b@@kb@x
  \ifvoid\p@p@nf@b@x\ifvoid\b@@nf@b@x
  \p@tb@x\v@lb@x\c@mm@ \ifvoid\y@@rb@x\else\sp@c@(\@nb@x\y@@rb@x)\fi
  \p@tb@x\is@b@x\c@mm@ \p@tb@x\p@g@b@x\c@mm@ \p@tb@x\inf@b@x\c@mm@
  \else\p@tb@x \b@@nf@b@x\c@mm@ \p@tb@x\v@lb@x\c@mm@ \p@tb@x\is@b@x\sp@c@
  \ifvoid\ed@b@x\else\sp@c@(\@nb@x\ed@b@x,\space\ifed@ ed.\else eds.\fi)\fi
  \p@tb@x\p@blb@x\c@mm@ \p@tb@x\p@bl@db@x\c@mm@ \p@tb@x\y@@rb@x\c@mm@
  \p@tb@x\p@g@b@x{\c@mm@\ifp@g@ p\p@@nt\else pp\p@@nt\fi}%
  \p@tb@x\inf@b@x\c@mm@\fi
  \else \p@tb@x\p@p@nf@b@x\c@mm@ \p@tb@x\v@lb@x\c@mm@
  \ifvoid\y@@rb@x\else\sp@c@(\@nb@x\y@@rb@x)\fi
  \p@tb@x\is@b@x\c@mm@ \p@tb@x\p@g@b@x\c@mm@ \p@tb@x\inf@b@x\c@mm@\fi
  \else \p@tb@x\b@@kb@x\c@mm@
  \p@tb@x\b@@nf@b@x\c@mm@ \p@tb@x\p@blb@x\c@mm@
  \p@tb@x\p@bl@db@x\c@mm@ \p@tb@x\y@@rb@x\c@mm@
  \ifvoid\p@g@b@x\else\c@mm@\@nb@x\p@g@b@x p\fi \p@tb@x\inf@b@x\c@mm@ \fi
  \else \c@mm@\@nb@x\p@p@rb@x\ic@\p@tb@x\p@p@nf@b@x\c@mm@
  \p@tb@x\v@lb@x\sp@c@ \ifvoid\y@@rb@x\else\sp@c@(\@nb@x\y@@rb@x)\fi
  \p@tb@x\is@b@x\c@mm@ \p@tb@x\p@g@b@x\c@mm@\p@tb@x\inf@b@x\c@mm@\fi
  \else \p@tb@x\p@p@rb@x\c@mm@\ic@\p@tb@x\p@p@nf@b@x\c@mm@
  \c@mm@\@nb@x\inb@@kb@x \p@tb@x\b@@nf@b@x\c@mm@ \p@tb@x\v@lb@x\sp@c@
  \p@tb@x\is@b@x\sp@c@
  \ifvoid\ed@b@x\else\sp@c@(\@nb@x\ed@b@x,\space\ifed@ ed.\else eds.\fi)\fi
  \p@tb@x\p@blb@x\c@mm@ \p@tb@x\p@bl@db@x\c@mm@ \p@tb@x\y@@rb@x\c@mm@
  \p@tb@x\p@g@b@x{\c@mm@\ifp@g@ p\p@@nt\else pp\p@@nt\fi}%
  \p@tb@x\inf@b@x\c@mm@\fi
  \else\p@tb@x\p@p@rb@x\c@mm@\ic@\p@tb@x\p@p@nf@b@x\c@mm@\p@tb@x\j@@rb@x\c@mm@
  \p@tb@x\v@lb@x\sp@c@ \ifvoid\y@@rb@x\else\sp@c@(\@nb@x\y@@rb@x)\fi
  \p@tb@x\is@b@x\c@mm@ \p@tb@x\p@g@b@x\c@mm@ \p@tb@x\inf@b@x\c@mm@ \fi}}
\def\m@r@f#1#2{\endr@f@\ifvoid\p@nct@\else\R@fb@x\hbox{\unhbox\r@f@b@x
 #1\unhbox\p@nct@\penalty-200\enskip#2}\fi\egroup\b@g@nr@f}
\def\endref{\endr@f@\ifvoid\p@nct@\else\R@fb@x\hbox{\unhbox\r@f@b@x.}\fi
 \parindent\r@f@nd
 \r@fit@m{\ifvoid\k@yb@x\else\if@fn@m\r@fn@{\unhbox\k@yb@x}\else
 \if@fl@b\unhbox\k@yb@x\else\r@fk@y{\unhbox\k@yb@x}\fi\fi\fi}\unhbox\r@f@b@x
 \endgraf\egroup\endgroup}
\def\moreref{\m@r@f;\empty}
\def\transl{\m@r@f;{\unskip\space
 {\sl English translation\ic@}:\penalty-66 \space}}
\def\endRefs{\endgraf\goodbreak\endgroup}

\hyphenation{acad-e-my acad-e-mies af-ter-thought anom-aly anom-alies
an-ti-deriv-a-tive an-tin-o-my an-tin-o-mies apoth-e-o-ses
apoth-e-o-sis ap-pen-dix ar-che-typ-al as-sign-a-ble as-sist-ant-ship
as-ymp-tot-ic asyn-chro-nous at-trib-uted at-trib-ut-able bank-rupt
bank-rupt-cy bi-dif-fer-en-tial blue-print busier busiest
cat-a-stroph-ic cat-a-stroph-i-cally con-gress cross-hatched data-base
de-fin-i-tive de-riv-a-tive dis-trib-ute dri-ver dri-vers eco-nom-ics
econ-o-mist elit-ist equi-vari-ant ex-quis-ite ex-tra-or-di-nary
flow-chart for-mi-da-ble forth-right friv-o-lous ge-o-des-ic
ge-o-det-ic geo-met-ric griev-ance griev-ous griev-ous-ly
hexa-dec-i-mal ho-lo-no-my ho-mo-thetic ideals idio-syn-crasy
in-fin-ite-ly in-fin-i-tes-i-mal ir-rev-o-ca-ble key-stroke
lam-en-ta-ble light-weight mal-a-prop-ism man-u-script mar-gin-al
meta-bol-ic me-tab-o-lism meta-lan-guage me-trop-o-lis
met-ro-pol-i-tan mi-nut-est mol-e-cule mono-chrome mono-pole
mo-nop-oly mono-spline mo-not-o-nous mul-ti-fac-eted mul-ti-plic-able
non-euclid-ean non-iso-mor-phic non-smooth par-a-digm par-a-bol-ic
pa-rab-o-loid pa-ram-e-trize para-mount pen-ta-gon phe-nom-e-non
post-script pre-am-ble pro-ce-dur-al pro-hib-i-tive pro-hib-i-tive-ly
pseu-do-dif-fer-en-tial pseu-do-fi-nite pseu-do-nym qua-drat-ic
quad-ra-ture qua-si-smooth qua-si-sta-tion-ary qua-si-tri-an-gu-lar
quin-tes-sence quin-tes-sen-tial re-arrange-ment rec-tan-gle
ret-ri-bu-tion retro-fit retro-fit-ted right-eous right-eous-ness
ro-bot ro-bot-ics sched-ul-ing se-mes-ter semi-def-i-nite
semi-ho-mo-thet-ic set-up se-vere-ly side-step sov-er-eign spe-cious
spher-oid spher-oid-al star-tling star-tling-ly sta-tis-tics
sto-chas-tic straight-est strange-ness strat-a-gem strong-hold
sum-ma-ble symp-to-matic syn-chro-nous topo-graph-i-cal tra-vers-a-ble
tra-ver-sal tra-ver-sals treach-ery turn-around un-at-tached
un-err-ing-ly white-space wide-spread wing-spread wretch-ed
wretch-ed-ly Brown-ian Eng-lish Euler-ian Feb-ru-ary Gauss-ian
Grothen-dieck Hamil-ton-ian Her-mit-ian Jan-u-ary Japan-ese Kor-te-weg
Le-gendre Lip-schitz Lip-schitz-ian Mar-kov-ian Noe-ther-ian
No-vem-ber Rie-mann-ian Schwarz-schild Sep-tem-ber}

\let\nopagenumber\p@gen@false \let\putpagenumber\p@gen@true

\outer\def\myRefs{\myR@fs} \r@st@re\proclaim
\def\bye{\par\vfill\supereject\cl@selbl\cl@secd\b@e} \r@endd@\b@e
 \let\Key\key \def\endpro{\par\endproclaim}
\let\d@c@\document \def\document{\d@c@\tenpoint}

\newtoks\@@tp@t \@@tp@t\output
\output=\@ft@{\let\{\noexpand\the\@@tp@t}
\let\{\relax

\newif\ifVersion \Versiontrue
\def\p@n@l#1{\ifnum#1=\z@\else\penalty#1\relax\fi}

\def\s@ct#1#2{\ifVersion
 \skip@\lastskip\ifdim\skip@<1.5\bls\vskip-\skip@\p@n@l{-200}\vsk.5>%
 \p@n@l{-200}\vsk.5>\p@n@l{-200}\vsk.5>\p@n@l{-200}\vsk-1.5>\else
 \p@n@l{-200}\fi\ifdim\skip@<.9\bls\vsk.9>\else
 \ifdim\skip@<1.5\bls\vskip\skip@\fi\fi
 \vtop{\twelvepoint\raggedright\s@cf@nt\vp1\vsk->\vskip.16ex
 \s@twd@\parindent{#1}%
 \ifdim\parindent>\z@\adv\parindent.5em\fi\hang\textindent{#1}#2\strut}
 \else
 \p@sk@p{-200}{.8\bls}\vtop{\s@cf@nt\s@twd@\parindent{#1}%
 \ifdim\parindent>\z@\adv\parindent.5em\fi\hang\textindent{#1}#2\strut}\fi
 \nointerlineskip\nobreak\vtop{\strut}\nobreak\vskip-.6\bls\nobreak}

\def\s@bs@ct#1#2{\ifVersion
 \skip@\lastskip\ifdim\skip@<1.5\bls\vskip-\skip@\p@n@l{-200}\vsk.5>%
 \p@n@l{-200}\vsk.5>\p@n@l{-200}\vsk.5>\p@n@l{-200}\vsk-1.5>\else
 \p@n@l{-200}\fi\ifdim\skip@<.9\bls\vsk.9>\else
 \ifdim\skip@<1.5\bls\vskip\skip@\fi\fi
 \vtop{\elevenpoint\raggedright\s@bf@nt\vp1\vsk->\vskip.16ex%
 \s@twd@\parindent{#1}\ifdim\parindent>\z@\adv\parindent.5em\fi
 \hang\textindent{#1}#2\strut}
 \else
 \p@sk@p{-200}{.6\bls}\vtop{\s@bf@nt\s@twd@\parindent{#1}%
 \ifdim\parindent>\z@\adv\parindent.5em\fi\hang\textindent{#1}#2\strut}\fi
 \nointerlineskip\nobreak\vtop{\strut}\nobreak\vskip-.8\bls\nobreak}

\def\gadv{\global\adv} \def\gad#1{\gadv#1\@ne} \def\gadneg#1{\gadv#1-\@ne}

\newcount\t@@n \t@@n=10 \newbox\testbox

\newcount\Sno \newcount\Lno \newcount\Fno

\def\pr@cl#1{\r@st@re\pr@c@\pr@c@{#1}\global\let\pr@c@\relax}

\def\l@L#1{\l@bel{#1}L} \def\l@F#1{\l@bel{#1}F} \def\<#1>{\l@b@l{#1}F}

\def\tagg#1{\tag"\rlap{\rm(#1)}\kern.01\p@"}
\def\Tag#1{\tag{\l@F{#1}}} \def\Tagg#1{\tagg{\l@F{#1}}}

\def\xspace{\kern.34em}

\def\Th#1{\pr@cl{Theorem\xspace\l@L{#1}}\ignore}
\def\Lm#1{\pr@cl{Lemma\xspace\l@L{#1}}\ignore}
\def\Cr#1{\pr@cl{Corollary\xspace\l@L{#1}}\ignore}
\def\Df#1{\pr@cl{Definition\xspace\l@L{#1}}\ignore}
\def\Cj#1{\pr@cl{Conjecture\xspace\l@L{#1}}\ignore}
\def\Prop#1{\pr@cl{Proposition\xspace\l@L{#1}}\ignore}
\def\Rem{\demo{\sl Remark}} 
\def\Pf#1.{\demo{Proof #1}} \def\epf{\qed\enddemo}

\def\Proof#1.{\demo{\let\{\relax Proof #1}\def\t@st@{#1}%
 \ifx\t@st@\empty\else\xdef\@@wr##1##2##3##4{##1{##2##3{\the\cdn@}{##4}}}%
 \wr@@c@{\the\cdn@}{Proof #1}\@@wr\wr@@c\string\subcd{\the\pageno}\fi\ignore}

\def\Ap@x{Appendix}
\def\Appendix{\Sno=64 \t@@n\@ne \wr@@c{\string\Appencd}
 \def\sf@rm{\char\the\Sno} \def\sf@rm@{\Ap@x\space\sf@rm} \def\sf@rm@@{\Ap@x}
 \def\s@ct@n##1##2{\s@ct\empty{\setbox\z@\hbox{##1}\ifdim\wd\z@=\z@
 \if##2*\sf@rm@@\else\if##2.\sf@rm@@.\else##2\fi\fi\else
 \if##2*\sf@rm@\else\if##2.\sf@rm@.\else\sf@rm@.\enspace##2\fi\fi\fi}}}
\def\Appcd#1#2#3{\gad\Cdentry\global\cdentry\z@\def\Ap@@{\hglue-\l@ftcd\Ap@x}
 \ifx\@ppl@ne\empty\def\l@@b{\@fwd@@{#1}{\space#1}{}}
 \if*#2\entcd{}{\Ap@@\l@@b}{#3}\else\if.#2\entcd{}{\Ap@@\l@@b.}{#3}\else
 \entcd{}{\Ap@@\l@@b.\enspace#2}{#3}\fi\fi\else
 \def\l@@b{\@fwd@@{#1}{\c@l@b{#1}}{}}\if*#2\entcd{\l@@b}{\Ap@x}{#3}\else
 \if.#2\entcd{\l@@b}{\Ap@x.}{#3}\else\entcd{\l@@b}{#2}{#3}\fi\fi\fi}

\let\s@ct@n\s@ct
\def\s@ct@@[#1]#2{\@ft@\xdef\csname @#1@S@\endcsname{\sf@rm}\wr@@x{}%
 \wr@@x{\string\labeldef{S}\space{\?#1@S?}\space{#1}}%
 {%\let\{\relax
 \s@ct@n{\sf@rm@}{#2}}\wr@@c{\string\Entcd{\?#1@S?}{#2}{\the\pageno}}}
\def\s@ct@#1{\wr@@x{}{%\let\{\relax
 \s@ct@n{\sf@rm@}{#1}}\wr@@c{\string\Entcd{\sf@rm}{#1}{\the\pageno}}}
\def\s@ct@e[#1]#2{\@ft@\xdef\csname @#1@S@\endcsname{\sf@rm}\wr@@x{}%
 \wr@@x{\string\labeldef{S}\space{\?#1@S?}\space{#1}}%
 {%\let\{\relax
 \s@ct@n\empty{#2}}\wr@@c{\string\Entcd{}{#2}{\the\pageno}}}
\def\s@cte#1{\wr@@x{}{%\let\{\relax
 \s@ct@n\empty{#1}}\wr@@c{\string\Entcd{}{#1}{\the\pageno}}}
\def\theSno#1#2{\dff\?#1@S?{#2}%
 \wr@@x{\string\labeldef{S}\space{#2}\space{#1}}\fi}

\newif\ifd@bn@\d@bn@true
\def\Section{\gad\Sno\ifd@bn@\Fno\z@\Lno\z@\fi\@fn@xt[\s@ct@@\s@ct@}
\def\section{\gad\Sno\ifd@bn@\Fno\z@\Lno\z@\fi\@fn@xt[\s@ct@e\s@cte}
\let\Sect\Section \let\sect\section
\def\subsection{\@fn@xt*\subs@ct@\subs@ct}
\def\subs@ct#1{{\s@bs@ct\empty{#1}}\wr@@c{\string\subcd{#1}{\the\pageno}}}
\def\subs@ct@*#1{\vsk->\nobreak
 {\s@bs@ct\empty{#1}}\wr@@c{\string\subcd{#1}{\the\pageno}}}

\def\l@b@l#1#2{\def\n@@{\csname #2no\endcsname}%
 \if*#1\gad\n@@ \@ft@\xdef\csname @#1@#2@\endcsname{\l@f@rm}\else\def\t@st{#1}%
 \ifx\t@st\empty\gad\n@@ \@ft@\xdef\csname @#1@#2@\endcsname{\l@f@rm}%
 \else\@ft@\ifx\csname @#1@#2@mark\endcsname\relax\gad\n@@
 \@ft@\xdef\csname @#1@#2@\endcsname{\l@f@rm}%
 \@ft@\gdef\csname @#1@#2@mark\endcsname{}%
 \wr@@x{\string\labeldef{#2}\space{\?#1@#2?}\space\ifnum\n@@<10 \space\fi{#1}}%
 \fi\fi\fi}
\def\labeldef#1#2#3{\dff\?#3@#1?{#2}}
\def\Labeldef#1#2#3{\dff\?#3@#1?{#2}\@ft@\gdef\csname @#3@#1@mark\endcsname{}}

\def\l@bel#1#2{\l@b@l{#1}{#2}\?#1@#2?}

\newcount\c@cite
\def\?#1?{\csname @#1@\endcsname}
\def\[{\@fn@xt:\c@t@sect\c@t@}
\def\c@t@#1]{{\c@cite\z@\@fwd@@{\?#1@L?}{\adv\c@cite1}{}%
 \@fwd@@{\?#1@F?}{\adv\c@cite1}{}\@fwd@@{\?#1?}{\adv\c@cite1}{}%
 \relax\ifnum\c@cite=\z@{\bf ???}\wrs@x{No label [#1]}\else
 \ifnum\c@cite=1\let\@@PS\relax\let\@@@\relax\else\let\@@PS\underbar
 \def\@@@{{\rm<}}\fi\@@PS{\?#1?\@@@\?#1@L?\@@@\?#1@F?}\fi}}
\def\(#1){{\rm(\c@t@#1])}}
\def\c@t@s@ct#1{\@fwd@@{\?#1@S?}{\?#1@S?\relax}%
 {{\bf ???}\wrs@x{No section label {#1}}}}
\def\c@t@sect:#1]{\c@t@s@ct{#1}} \let\SNo\c@t@s@ct

\newdimen\l@ftcd \newdimen\r@ghtcd \let\nlc\relax
\newcount\Cdentry \newcount\cdentry \let\prentcd\relax \let\postentcd\relax

\def\d@tt@d{\leaders\hbox to 1em{\kern.1em.\hfil}\hfill}
\def\entcd#1#2#3{\gad\cdentry\prentcd\item{\l@bcdf@nt#1}{\entcdf@nt#2}\alb
 \kern.9em\hbox{}\kern-.9em\d@tt@d\kern-.36em{\p@g@cdf@nt#3}\kern-\r@ghtcd
 \hbox{}\postentcd\par}
\def\Entcd#1#2#3{\gad\Cdentry\global\cdentry\z@
 \def\l@@b{\@fwd@@{#1}{\c@l@b{#1}}{}}\vsk.2>\entcd{\l@@b}{#2}{#3}}
\def\subcd#1#2{{\adv\leftskip.333em\entcd{}{\s@bcdf@nt#1}{#2}}}
\def\Refcd#1#2{\def\t@@st{#1}\ifx\t@@st\empty\ifx\r@fl@ne\empty\relax\else
 \R@fcd{\r@fl@ne}{#2}\fi\else\R@fcd{#1}{#2}\fi}
\def\R@fcd#1#2{\sk@@p{.6\bls}\entcd{}{\hglue-\l@ftcd\R@fcdf@nt #1}{#2}}
\def\Refline{\def\r@fl@ne} \def\Refempty{\let\r@fl@ne\empty}
\def\Appencd{\par\adv\leftskip-\l@ftcd\adv\rightskip-\r@ghtcd\@ppl@ne
 \adv\leftskip\l@ftcd\adv\rightskip\r@ghtcd\let\Entcd\Appcd}
\def\appline{\def\@ppl@ne} \def\Appempty{\let\@ppl@ne\empty}
\def\Appline#1{\def\@ppl@ne{\s@bs@ct{}{#1}}}
\def\Leftcd#1{\adv\leftskip-\l@ftcd\s@twd@\l@ftcd{\c@l@b{#1}\enspace}
 \adv\leftskip\l@ftcd}
\def\Rightcd#1{\adv\rightskip-\r@ghtcd\s@twd@\r@ghtcd{#1\enspace}
 \adv\rightskip\r@ghtcd}
\def\C@nt{Contents} \def\Ap@s{Appendices} \def\R@fcs{References}
\def\contents{\@fn@xt*\cont@@\cont@}
\def\cont@{\@fn@xt[\cnt@{\cnt@[\C@nt]}}
\def\cont@@*{\@fn@xt[\cnt@@{\cnt@@[\C@nt]}}
\def\cnt@[#1]{\c@nt@{M}{#1}{44}{\s@bs@ct{}{\@ppl@f@nt\Ap@s}}}
\def\cnt@@[#1]{\c@nt@{M}{#1}{44}{}}
\def\endco{\par\penalty-500\vsk>\vskip\z@\endgroup}
\def\readcd{\@np@t{\jobname.cd}}
\def\Cde{\@fn@xt*\Cde@@\Cde@}
\def\Cde@{\@fn@xt[\Cd@{\Cd@[\C@nt]}}
\def\Cde@@*{\@fn@xt[\Cd@@{\Cd@@[\C@nt]}}
\def\Cd@[#1]{\cnt@[#1]\readcd\endco}
\def\Cd@@[#1]{\cnt@@[#1]\readcd\endco}
\def\contlabeldef{\def\c@l@b}

\long\def\c@nt@#1#2#3#4{\s@twd@\l@ftcd{\c@l@b{#1}\enspace}
 \s@twd@\r@ghtcd{#3\enspace}\adv\r@ghtcd1.333em
 \def\@ppl@ne{#4}\def\r@fl@ne{\R@fcs}\s@ct{}{#2}\B@gr@\parindent\z@\let\nlc\nl
 \let\nl\relax\parskip.2\bls\adv\leftskip\l@ftcd\adv\rightskip\r@ghtcd}

\def\writecd{\immediate\openout\@@cd\jobname.cd \def\wr@@c{\write\@@cd}
 \def\cl@secd{\immediate\write\@@cd{\string\endinput}\immediate\closeout\@@cd}
 \def\closecd{\cl@secd\global\let\cl@secd\relax}}
\let\cl@secd\relax \def\wr@@c#1{} \let\closecd\relax

\def\dff{\@ft@\d@f} \def\d@f{\@ft@\def}
\def\edff{\@ft@\ed@f} \def\ed@f{\@ft@\edef}
\def\gdff{\@ft@\gd@f} \def\gd@f{\@ft@\gdef}
\def\defi#1#2{\def#1{#2}\wr@@x{\string\def\string#1{#2}}}

\def\qed{\hbox{}\nobreak\hfill\nobreak{\m@th$\,\square$}}
\def\back#1 {\strut\kern-.33em #1\enspace\ignore} %% !!! a space after #1 !!!
\def\Text#1{\crcr\noalign{\alb\vsk>\normalbaselines\vsk->\vbox{\nt #1\strut}%
 \nobreak\nointerlineskip\vbox{\strut}\nobreak\vsk->\nobreak}}

\def\hcor#1{\advance\hoffset by #1}
\def\vcor#1{\advance\voffset by #1}
\let\bls\baselineskip \let\ignore\ignorespaces
\ifx\ic@\undefined \let\ic@\/\fi
\def\vsk#1>{\vskip#1\bls} \let\adv\advance
\def\vv#1>{\vadjust{\vsk#1>}\ignore}
\def\vvn#1>{\vadjust{\nobreak\vsk#1>\nobreak}\ignore}
\def\vvv#1>{\vskip\z@\vsk#1>\nt\ignore}
\def\vvgood{\vadjust{\penalty-500}}
\def\nngood{\noalign{\penalty-500}}

\def\Goodbreak{\par\penalty-\@m}
\def\wgood#1>{\vv#1>\vvgood\vv-#1>}
\def\wwgood#1:#2>{\vv#1>\vvgood\vv#2>}
\def\mmgood#1:#2>{\cnn#1>\nngood\cnn#2>}
\def\goodsk#1:#2>{\vsk#1>\goodbreak\vsk#2>\vsk0>}
\def\ragood{\vadjust{\vskip\z@ plus 12pt}\vvgood}

\def\Par{\vsk.5>} \def\setparindent{\edef\Parindent{\the\parindent}}
\def\Type{\vsk.5>\bgroup\parindent\z@\tt\rightskip\z@ plus1em minus1em%
 \spaceskip.3333em \xspaceskip.5em\relax}
\def\endType{\vsk.5>\egroup\nt} 

  \let\dollar\$ \let\ampersand\&
\let\sss\scriptscriptstyle  
\let\vp\vphantom \let\hp\hphantom \let\nt\noindent
\let\cline\centerline \let\lline\leftline \let\rline\rightline
\def\nn#1>{\noalign{\vskip#1\p@@}} \def\NN#1>{\openup#1\p@@}
\def\cnn#1>{\noalign{\vsk#1>}}
 
\let\Lim\lim \def\lim{\Lim\limits} \let\Sum\sum \def\sum{\Sum\limits}
\def\Plus{\bigoplus\limits} 
\let\Prod\prod \def\prod{\Prod\limits} \let\Int\int \def\int{\Int\limits}

\def\tsum{\mathop{\tsize\Sum}\limits} 
\def\tprod{\mathop{\tsize\Prod}\limits} \def\&{.\kern.1em}
\def\nl{\leavevmode\hfill\break} \def\~{\leavevmode\@fn@xt~\m@n@s\@md@@sh}
\def\@md@@sh{\@fn@xt-\d@@sh\@md@sh} \def\@md@sh{\raise.13ex\hbox{--}}
\def\m@n@s~{\raise.15ex\mbox{-}} \def\d@@sh-{\raise.15ex\hbox{-}}

\let\procent\% \def\%#1{\ifmmode\mathop{#1}\limits\else\procent#1\fi}
\let\@ml@t\" \def\"#1{\ifmmode ^{(#1)}\else\@ml@t#1\fi}
\let\@c@t@\' \def\'#1{\ifmmode _{(#1)}\else\@c@t@#1\fi}
\let\colon\: \def\:{^{\vp{\topsmash|}}} 

\let\texspace\ \def\ {\ifmmode\alb\fi\texspace} \def\.{\d@t\ignore}
%%% Do not remove a space after \ !!!

\newif\ifNewskips

\def\Newskips{\global\Newskipstrue
 \gdef\>{\RIfM@\mskip.666667\thinmuskip\relax\else\kern.111111em\fi}
 \gdef\}{\RIfM@\mskip-.666667\thinmuskip\relax\else\kern-.111111em\fi}
 \gdef\){\RIfM@\mskip.333333\thinmuskip\relax\else\kern.0555556em\fi}
 \gdef\]{\RIfM@\mskip-.333333\thinmuskip\relax\else\kern-.0555556em\fi}}
\def\d@t{\ifNewskips.\hskip.3em\else\def\d@t{.\ }\fi} \def\.{\d@t\ignore}
\Newskips

\let\n@wp@ge\newpage \def\newpage{\endgraf\n@wp@ge}
\let\=\m@th \def\mbox#1{\hbox{\m@th$#1$}}
\def\mtext#1{\text{\m@th$#1$}} \def\^#1{\text{\m@th#1}}
\def\Line#1{\kern-.5\hsize\line{\m@th$\dsize#1$}\kern-.5\hsize}
\def\Lline#1{\kern-.5\hsize\lline{\m@th$\dsize#1$}\kern-.5\hsize}
\def\Cline#1{\kern-.5\hsize\cline{\m@th$\dsize#1$}\kern-.5\hsize}
\def\Rline#1{\kern-.5\hsize\rline{\m@th$\dsize#1$}\kern-.5\hsize}

\def\Ll@p#1{\llap{\m@th$#1$}} \def\Rl@p#1{\rlap{\m@th$#1$}}
 \def\Cl@p#1{\llap{\m@th$#1$\hss}}
\def\Llap#1{\mathchoice{\Ll@p{\dsize#1}}{\Ll@p{\tsize#1}}{\Ll@p{\ssize#1}}%
 {\Ll@p{\sss#1}}}
\def\Clap#1{\mathchoice{\Cl@p{\dsize#1}}{\Cl@p{\tsize#1}}{\Cl@p{\ssize#1}}%
 {\Cl@p{\sss#1}}}
\def\Rlap#1{\mathchoice{\Rl@p{\dsize#1}}{\Rl@p{\tsize#1}}{\Rl@p{\ssize#1}}%
 {\Rl@p{\sss#1}}}
 
\def\LRtph#1#2{\setbox\z@\hbox{#1}\dimen\z@\wd\z@\hbox{\hbox to\dimen\z@{#2}}}
\def\LRph#1#2{\LRtph{\m@th$#1$}{\m@th$#2$}}

\def\LLph#1#2{\LRph{#1}{\hss#2}} 

\def\Lph#1#2{\mathchoice{\LLph{\dsize#1}{\dsize#2}}{\LLph{\tsize#1}{\tsize#2}}
 {\LLph{\ssize#1}{\ssize#2}}{\LLph{\sss#1}{\sss#2}}}

\def\Lto#1{\setbox\z@\mbox{\tsize{#1}}%
 \mathrel{\mathop{\hbox to\wd\z@{\rightarrowfill}}\limits#1}}
\def\Lgets#1{\setbox\z@\mbox{\tsize{#1}}%
 \mathrel{\mathop{\hbox to\wd\z@{\leftarrowfill}}\limits#1}}
\def\vpb#1{{\vp{\big(}}^{\]#1}} \def\vpp#1{{\vp{\big]}}_{#1}}

\let\alb\allowbreak 
\def\ald{\noalign{\alb}} 

\let\o\circ \let\x\times \let\ox\otimes 
\let\sub\subset  \let\tabs\+
\let\le\leqslant \let\ge\geqslant
\let\der\partial \let\8\infty \let\*\star
\let\bra\langle \let\ket\rangle
 
\let\map\mapsto  \let\hto\hookrightarrow
 
\let\Vert\parallel  \def\nin{\not\in}

\let\lb\lbrace \let\rb\rbrace

\def\lsym#1{#1\alb\ldots\relax#1\alb}
\def\lc{\lsym,}   \def\lox{\lsym\ox}

\def\End{\mathop{\roman{End}\>}\nolimits}

\def\id{\roman{id}}  
 \def\1{^{-1}} \let\underscore\_ \def\_#1{_{\Rlap{#1}}}
\def\vst#1{{\lower1.9\p@@\mbox{\bigr|_{\raise.5\p@@\mbox{\ssize#1}}}}}
\def\vrp#1:#2>{{\vrule height#1 depth#2 width\z@}}
\def\vru#1>{\vrp#1:\z@>} \def\vrd#1>{\vrp\z@:#1>}
\def\qqq{\qquad\quad} 
\def\sscr#1{\raise.3ex\mbox{\sss#1}} \def\@@PS{\bold{OOPS!!!}}

\def\intcl{\mathop
 {\Rlap{\raise.3ex\mbox{\kern.12em\curvearrowleft}}\int}\limits}
\def\intcr{\mathop
 {\Rlap{\raise.3ex\mbox{\kern.24em\curvearrowright}}\int}\limits}

\def\pms{\raise.25ex\mbox{\ssize\pm}\>}
\def\mps{\raise.25ex\mbox{\ssize\mp}\>}
\def\pss{{\sscr+}} \def\mss{{\sscr-}}

\let\al\alpha

 \let\Gm\Gamma 
\let\dl\delta  
 \let\eps\varepsilon \let\epsilon\eps

\let\ka\kappa
\let\la\lambda

 \let\phi\varphi

 \let\Om\Omega 

\def\C{\Bbb C}

\def\Z{\Bbb Z}

\def\DD{\Bbb D}

\def\Zp{\Z_{\ge 0}}

\def\difl/{differential} \def\dif/{difference}
\def\cf.{cf.\ \ignore} \def\Cf.{Cf.\ \ignore}
\def\egv/{eigenvector} \def\eva/{eigenvalue} \def\eq/{equation}
\def\lhs/{the left hand side} \def\rhs/{the right hand side}
\def\Lhs/{The left hand side} \def\Rhs/{The right hand side}
\def\gby/{generated by} \def\wrt/{with respect to} \def\st/{such that}
\def\resp/{respectively} \def\off/{offdiagonal} \def\wt/{weight}
\def\pol/{polynomial} \def\rat/{rational} \def\tri/{trigonometric}
\def\fn/{function} \def\var/{variable} \def\raf/{\rat/ \fn/}
\def\inv/{invariant} \def\hol/{holomorphic} \def\hof/{\hol/ \fn/}
\def\mer/{meromorphic} \def\mef/{\mer/ \fn/} \def\mult/{multiplicity}
\def\sym/{symmetric} \def\perm/{permutation}
\def\rep/{representation} \def\irr/{irreducible} \def\irrep/{\irr/ \rep/}
\def\hom/{homomorphism} \def\aut/{automorphism} \def\iso/{isomorphism}
\def\lex/{lexicographical} \def\as/{asymptotic} \def\asex/{\as/ expansion}
\def\ndeg/{nondegenerate} \def\neib/{neighbourhood} \def\deq/{\dif/ \eq/}
\def\hw/{highest \wt/} \def\gv/{generating vector} \def\eqv/{equivalent}
\def\msd/{method of steepest descend} \def\pd/{pairwise distinct}
\def\wlg/{without loss of generality} \def\Wlg/{Without loss of generality}
\def\onedim/{one-dim\-en\-sion\-al} \def\fd/{fin\-ite-dim\-en\-sion\-al}
\def\qcl/{quasiclassical} \def\hwv/{\hw/ vector}
\def\hgeom/{hyper\-geo\-met\-ric} \def\hint/{\hgeom/ integral}
\def\hwm/{\hw/ module} \def\emod/{evaluation module} \def\Vmod/{Verma module}
\def\symg/{\sym/ group} \def\sol/{solution} \def\eval/{evaluation}
\def\anf/{analytic \fn/} \def\anco/{analytic continuation}
\def\qg/{quantum group} \def\qaff/{quantum affine algebra}

\hyphenation{ortho-gon-al}

\def\Rm/{\^{$R$-}matrix} \def\Rms/{\^{$R$-}matrices} \def\YB/{Yang-Baxter \eq/}
\def\Ba/{Bethe ansatz} \def\Bv/{Bethe vector} \def\Bae/{\Ba/ \eq/}
\def\KZv/{Knizh\-nik-Zamo\-lod\-chi\-kov} \def\KZvB/{\KZv/-Bernard}
\def\KZ/{{\sl KZ\/}} \def\qKZ/{{\sl qKZ\/}}
\def\KZB/{{\sl KZB\/}} \def\qKZB/{{\sl qKZB\/}}
\def\qKZo/{\qKZ/ operator} \def\qKZc/{\qKZ/ connection}
\def\KZe/{\KZ/ \eq/} \def\qKZe/{\qKZ/ \eq/} \def\qKZBe/{\qKZB/ \eq/}
\def\XXX/{{\sl XXX\/}} \def\XXZ/{{\sl XXZ\/}} \def\XYZ/{{\sl XYZ\/}}

\def\h@ph{\discretionary{}{}{-}} \def\$#1$-{\,\^{$#1$}\h@ph}

\def\TFT/{Research Insitute for Theoretical Physics}
\def\HY/{University of Helsinki} \def\AoF/{the Academy of Finland}
\def\CNRS/{Supported in part by MAE\~MICECO\~CNRS Fellowship}
\def\LPT/{Laboratoire de Physique Th\'eorique ENSLAPP}
\def\ENSLyon/{\'Ecole Normale Sup\'erieure de Lyon}
\def\LPTaddr/{46, All\'ee d'Italie, 69364 Lyon Cedex 07, France}
\def\enslapp/{URA 14\~36 du CNRS, associ\'ee \`a l'E.N.S.\ de Lyon,
au LAPP d'Annecy et \`a l'Universit\`e de Savoie}
\def\ensemail/{vtarasov\@ enslapp.ens-lyon.fr}
\def\DMS/{Department of Mathematics, Faculty of Science}
\def\DMO/{\DMS/, Osaka University}
\def\DMOaddr/{Toyonaka, Osaka 560, Japan}
\def\dmoemail/{vt\@ math.sci.osaka-u.ac.jp}
\def\MPI/{Max\)-Planck\)-Institut} \def\MPIM/{\MPI/ f\"ur Mathematik}
\def\MPIMaddr/{P\]\&O.\ Box 7280, D\~-\]53072 \,Bonn, Germany}
\def\mpimemail/{tarasov\@ mpim-bonn.mpg.de}
\def\SPb/{St\&Peters\-burg}
\def\home/{\SPb/ Branch of Steklov Mathematical Institute}
\def\homeaddr/{Fontanka 27, \SPb/ \,191011, Russia}
\def\homemail/{vt\@ pdmi.ras.ru}
\def\absence/{On leave of absence from \home/}
\def\support/{Supported in part by}
\def\UNC/{Department of Mathematics, University of North Carolina}
\def\ChH/{Chapel Hill} \def\UNCaddr/{\ChH/, NC 27599, USA}
\def\avemail/{anv\@ email.unc.edu}	%% {av\@ math.unc.edu}
\def\grant/{NSF grant DMS\~9501290}	%% Felder's grant no. 9400841
\def\Grant/{\support/ \grant/}

\def\Aomoto/{K\&Aomoto}
\def\Cher/{I\&Che\-red\-nik}
\def\Dri/{V\]\&G\&Drin\-feld}
\def\Fadd/{L\&D\&Fad\-deev}
\def\Feld/{G\&Felder}
\def\Fre/{I\&B\&Fren\-kel}
\def\Etingof/{P\]\&Etingof}
\def\Gustaf/{R\&A\&Gustafson}
\def\Izergin/{A\&G\&Izergin}
\def\Jimbo/{M\&Jimbo}
\def\Kazh/{D\&Kazhdan}
\def\Kor/{V\]\&E\&Kore\-pin}
\def\Kulish/{P\]\&P\]\&Ku\-lish}
\def\Lusz/{G\&Lusztig}
\def\Miwa/{T\]\&Miwa}
\def\MN/{M\&Naza\-rov}
\def\Reshet/{N\&Reshe\-ti\-khin} \def\Reshy/{N\&\]Yu\&Reshe\-ti\-khin}
\def\SchV/{V\]\&\]V\]\&Schecht\-man} \def\Sch/{V\]\&Schecht\-man}
\def\Skl/{E\&K\&Sklya\-nin}
\def\Smirnov/{F\]\&A\&Smir\-nov}
\def\Takh/{L\&A\&Takh\-tajan}
\def\VT/{V\]\&Ta\-ra\-sov} \def\VoT/{V\]\&O\&Ta\-ra\-sov}
\def\Varch/{A\&\]Var\-chenko} \def\Varn/{A\&N\&\]Var\-chenko}
\def\Zhel/{D\&P\]\&Zhe\-lo\-ben\-ko}

\def\AiA/{Al\-geb\-ra i Ana\-liz}
\def\DAN/{Do\-kla\-dy AN SSSR}
\def\FAA/{Funk\.Ana\-liz i ego pril.}
\def\Izv/{Iz\-ves\-tiya AN SSSR, ser\.Ma\-tem.}
\def\TMF/{Teo\-ret\.Ma\-tem\.Fi\-zi\-ka}
\def\UMN/{Uspehi Matem.\ Nauk}

\def\AMS/{Amer\.Math\.Society}
\def\AMSa/{AMS \publaddr Providence RI}
\def\AMST/{\AMS/ Transl.,\ Ser\&\)2}
\def\AMSTr/{\AMS/ Transl.,} \def\Ser2{Ser\&\)2}
\def\Astq/{Ast\'erisque}
\def\ContM/{Contemp\.Math.}
\def\CMP/{Comm\.Math\.Phys.}
\def\DMJ/{Duke\.Math\.J.}
\def\Inv/{Invent\.Math.} %% Inventiones Mathematicae
\def\IMRN/{Int\.Math\.Res.\ Notices}
\def\JMP/{J\.Math\.Phys.}
\def\JPA/{J\.Phys.\ A}
\def\JSM/{J\.Soviet Math.}
\def\LMJ/{Leningrad Math.\ J.}
\def\LpMJ/{\SPb/ Math.\ J.}
\def\LMP/{Lett\.Math\.Phys.}
\def\NMJ/{Nagoya Math\.J.}
\def\Nucl/{Nucl\.Phys.\ B}
\def\OJM/{Osaka J\.Math.}
\def\RIMS/{Publ\.RIMS, Kyoto Univ.}
\def\SIAM/{SIAM J\.Math\.Anal.}
\def\SMNS/{Selecta Math., New Series}
\def\TMP/{Theor\.Math\.Phys.}
\def\ZNS/{Zap\. nauch\. semin. LOMI}

\def\ASMP/{Advanced Series in Math.\ Phys.{}}

\def\Birk/{Birkh\"auser}
\def\CUP/{Cambridge University Press} \def\CUPa/{\CUP/ \publaddr Cambridge}
\def\Spri/{Springer\)-Verlag} \def\Spria/{\Spri/ \publaddr Berlin}
\def\WS/{World Scientific} \def\WSa/{\WS/ \publaddr Singapore}

\newbox\lefthbox \newbox\righthbox

\let\sectsep. \let\labelsep. \let\contsep. \let\labelspace\relax
\let\sectpre\relax \let\contpre\relax
\def\sf@rm{\the\Sno} \def\sf@rm@{\sectpre\sf@rm\sectsep}
\def\c@l@b#1{\contpre#1\contsep}
\def\l@f@rm{\ifd@bn@\sf@rm\labelsep\fi\labelspace\the\n@@}

\def\sectformdef{\def\sf@rm}

\let\DoubleNum\d@bn@true \let\SingleNum\d@bn@false

\def\NoNewNum{\let\writeldf\relax\def\l@b@l##1##2{\if*##1%
 \@ft@\xdef\csname @##1@##2@\endcsname{\mbox{*{*}*}}\fi}}
\def\NoNewTime{\def\todaydef##1{\def\today{##1}}
 \def\nowtimedef##1{\def\nowtime{##1}}}
\def\NoInput{\let\Input\input\let\writeldf\relax}
\def\Fixed{\NoNewTime\NoInput}

\newbox\dtlb@x
\def\DateTimeLabel{\global\setbox\dtlb@x\vbox to\z@{\ifMag\eightpoint\else
 \ninepoint\fi\sl\vss\rline\today\rline\nowtime}
 \global\headline{\hfil\box\dtlb@x}}

\def\sectfont#1{\def\s@cf@nt{#1}} \sectfont\bf
\def\subsectfont#1{\def\s@bf@nt{#1}} \subsectfont\it
\def\Entcdfont#1{\def\entcdf@nt{#1}} \Entcdfont\relax
\def\labelcdfont#1{\def\l@bcdf@nt{#1}} \labelcdfont\relax
\def\pagecdfont#1{\def\p@g@cdf@nt{#1}} \pagecdfont\relax
\def\subcdfont#1{\def\s@bcdf@nt{#1}} \subcdfont\it
\def\applefont#1{\def\@ppl@f@nt{#1}} \applefont\bf
\def\Refcdfont#1{\def\R@fcdf@nt{#1}} \Refcdfont\bf

\def\reffont#1{\def\r@ff@nt{#1}} \reffont\rm
\def\keyfont#1{\def\k@yf@nt{#1}} \keyfont\rm
\def\paperfont#1{\def\p@p@rf@nt{#1}} \paperfont\it
\def\bookfont#1{\def\b@@kf@nt{#1}} \bookfont\it
\def\volfont#1{\def\v@lf@nt{#1}} \volfont\bf
\def\issuefont#1{\def\iss@f@nt{#1}} \issuefont{no\p@@nt}

\def\adjustmid#1{\kern-#1\p@\alb\hskip#1\p@\relax}
\def\adjustend#1{\adjustnext{\kern-#1\p@\alb\hskip#1\p@}}

\newif\ifcd 

\tenpoint

\Fixed

\loadbold

\font@\Beufm=eufm10 scaled 1440
\font@\Eufm=eufm8 scaled 1440
\newfam\Frakfam \textfont\Frakfam\Eufm

\def\xspace{\kern.34em}

\def\avemail/{anv\@ email.unc.edu}

\def\Th#1{\pr@cl{Theorem\xspace\l@L{#1}}\ignore}
\def\Lm#1{\pr@cl{Lemma\xspace\l@L{#1}}\ignore}
\def\Cr#1{\pr@cl{Corollary\xspace\l@L{#1}}\ignore}
\def\Df#1{\pr@cl{Definition\xspace\l@L{#1}}\ignore}
\def\Cj#1{\pr@cl{Conjecture\xspace\l@L{#1}}\ignore}
\def\Prop#1{\pr@cl{Proposition\xspace\l@L{#1}}\ignore}
\def\Pf#1.{\demo{Proof #1}\bgroup\ignore}
\def\epf{\qed\par\egroup\enddemo}

\def\Vert{\ \,\big|\ \,}

\def\##1{^{[#1]}}

\def\db{\bold d}

\def\Pc{\Cal P}
\def\Zc{\Cal Z}

\def\gg{\frak g}
\def\hg{\frak h}
\def\ng{\frak n}

\def\gl{\frak{gl}}
\def\gsl{\frak{sl}}

\def\Dh{\widehat D}
\def\habla{\widehat\nabla}

\def\ev{\mathop{\slanted{ev}}\nolimits}

\def\ngp{\ng_{\sss+}} \def\ngm{\ng_{\sss-}} 

\def\Ck{\C^{\)k}} 
\def\glt{\gl_2} \def\glk{\gl_k} \def\gln{\gl_{\)n}}
\def\glkn{(\)\glk\>,\)\gln\))}
\def\Uqkn{\bigl(\)U_q(\glk)\>,\)U_q(\gln)\bigr)}
\def\slk{\gsl_{\)k}} \def\sln{\gsl_n}
\def\Uk{U(\glk)} 
\def\Yk{Y(\glk)} \def\Yn{Y(\gln)}
\def\Ukn{\bigl(U(\glk)\bigr)\vpb{\ox\)n}}

\def\Vb{V_\bullet}

\def\+#1{^{\sscr{\bra\]#1\]\ket}}}

\def\ak{a=1\lc k} \def\bk{b=1\lc k} \def\alk{\al_1\lc\al_{k-1}}
\def\inn{i=1\lc n} \def\jn{j=1\lc n}
\def\zn{z_1\lc z_n} \def\lak{\la_1\lc\la_k}
\def\xkn{x_{11}\lc x_{k1}\lc x_{1n}\)\lc x_{kn}}

\def\Dlk{D_{\la_1}\lc D_{\la_k}}
\def\Dhk{\Dh_{\la_1}\lc\Dh_{\la_k}}
\def\Qlk{Q_{\la_1}\lc Q_{\la_k}}
\def\Zzn{Z_{z_1}\lc Z_{z_n}}
\def\nablasn{\nabla_{\}z_1}\lc\nabla_{\}z_n}}
\def\hablasn{\habla_{\}z_1}\lc\habla_{\}z_n}}

\def\Omo{\Om^{\sscr{\)0}}} \def\Omp{\Om^\pss} \def\Omm{\Om^\mss}

\def\Pck{\Pc_k} \def\Pcm{\Pc_{\]m}} 
\def\Plm{P_{kn}[\)\la\>,\mu\)]}

\def\Zlm{\Zc_{kn}[\)\la\>,\mu\)]}

\def\Wla{W[\)\la\)]} 
\def\Vlt{V_{l_1}\+2\lox V_{l_n}\+2}
\def\Vmn{V\+n_{m_1}\ox V\+n_{m_2}} \def\Vmmn{V\+n\'{m_1+\)m_2-\)m,m,0\lc 0}}
\def\Vlm{(V_{l_1}\+k\lox V_{l_n}\+k)\)[\)m_1\lc m_k]}
\def\Vml{(V_{m_1}\+n\lox V_{m_k}\+n)\)[\)l_1\lc l_n]}
\def\Vlmt{(V_{l_1}\+2\lox V_{l_n}\+2)\)[\)m_1,m_2]}
\def\Vlmtm{\Vlmt\vpp{\)m}}
\def\Vmlt{(\Vmn)\)[\)l_1\lc l_n]} \def\Vmltm{\Vmlt\vpp{\)m}}
\def\VtWm{V\+2\'{|\la|-\)m\),\)m}\ox W\+2_m}
\def\Vbkn{\bigl(\)\Vb\+k\)\bigr)\vpb{\ox n}}
\def\Vbnk{\bigl(\)\Vb\+n\)\bigr)\vpb{\ox k}}

\def\PVW{P_{VW}\:}
\def\RVb{R_{\Vb\Vb\]}\:} \def\RVlm{R_{V_lV_m\]}\:}
\def\RUV{R_{UV}\:} \def\RUW{R_{UW}\:} \def\RVW{R_{VW}\:}

\def\DD/{{\sl DD\/}}
\def\qDD/{{\sl q\)DD\/}}

\def\ehom/{\eval/ \hom/} \def\part/{partition}

\def\gtmod/{\$\glt$-module}
\def\gkmod/{\$\glk$-module} \def\gnmod/{\$\gln$-module}
\def\skmod/{\$\slk$-module} \def\snmod/{\$\sln$-module}
\def\Ykmod/{\$\Yk$-module} \def\Ynmod/{\$\Yn$-module}

\let\goodbm\relax  \let\mmgood\relax \def\mngood{}
   
\let\mline\relax \def\vvm#1>{\ignore} \def\vvnm#1>{\ignore} \def\cnnm#1>{}
\def\cnnu#1>{} \def\vvu#1>{\ignore} \def\vvnu#1>{\ignore} 

\ifMag \ifUS   
  \let\vvu\vv \let\vvnu\vvn \let\cnnu\cnn  \else
 \let\goodbm\goodbreak  \let\mmgood\vvgood \let\cnnm\cnn
 \let\mline\nl \let\vvm\vv \let\vvnm\vvn \let\mngood\nngood \fi
 \let\goodbreak\relax  \let\vvgood\relax
  \def\nngood{}  \fi

\def\wwgood#1:#2>{\vv#1>\vvgood\vv#2>\vv0>}
\def\vskgood#1:#2>{\vsk#1>\goodbreak\vsk#2>\vsk0>}

\def\wwmgood#1:#2>{\ifMag\vv#1>\mmgood\vv#2>\vv0>\fi}
\def\vskmgood#1:#2>{\ifMag\vsk#1>\goodbm\vsk#2>\vsk0>\fi}
\def\vskm#1:#2>{\ifMag\vsk#1>\else\vsk#2>\fi}
\def\vvmm#1:#2>{\ifMag\vv#1>\else\vv#2>\fi}
\def\vvnn#1:#2>{\ifMag\vvn#1>\else\vvn#2>\fi}
\def\nnm#1:#2>{\ifMag\nn#1>\else\nn#2>\fi}
\def\kerm#1:#2>{\ifMag\kern#1em\else\kern#2em\fi}

\csname dual.def\endcsname

\labeldef{F} {1\labelsep \labelspace 1}  {Vb}

\labeldef{L} {2\labelsep \labelspace 1}  {ratKZ}
\labeldef{L} {2\labelsep \labelspace 2}  {trigKZ}
\labeldef{L} {2\labelsep \labelspace 3}  {KZDD}
\labeldef{L} {2\labelsep \labelspace 4}  {trigDD}

\labeldef{F} {3\labelsep \labelspace 1}  {Xa}
\labeldef{L} {3\labelsep \labelspace 1}  {KZqDD}

\labeldef{F} {4\labelsep \labelspace 1}  {Rinv}
\labeldef{F} {4\labelsep \labelspace 2}  {Rdef}
\labeldef{F} {4\labelsep \labelspace 3}  {inv}
\labeldef{F} {4\labelsep \labelspace 4}  {YB}
\labeldef{F} {4\labelsep \labelspace 5}  {Ki}
\labeldef{L} {4\labelsep \labelspace 1}  {qKZ}
\labeldef{L} {4\labelsep \labelspace 2}  {qKZDD}

\labeldef{F} {5\labelsep \labelspace 1}  {left}
\labeldef{F} {5\labelsep \labelspace 2}  {right}
\labeldef{F} {5\labelsep \labelspace 3}  {kiso}
\labeldef{F} {5\labelsep \labelspace 4}  {niso}
\labeldef{L} {5\labelsep \labelspace 1}  {Pkn}
\labeldef{L} {5\labelsep \labelspace 2}  {dual}
\labeldef{L} {5\labelsep \labelspace 3}  {Plm}
\labeldef{L} {5\labelsep \labelspace 4}  {basisk}
\labeldef{L} {5\labelsep \labelspace 5}  {basisn}
\labeldef{F} {5\labelsep \labelspace 5}  {vkx}
\labeldef{F} {5\labelsep \labelspace 6}  {vnx}
\labeldef{L} {5\labelsep \labelspace 6}  {PV}
\labeldef{L} {5\labelsep \labelspace 7}  {Ikn}
\labeldef{F} {5\labelsep \labelspace 7}  {C}
\labeldef{L} {5\labelsep \labelspace 8}  {dualKZDD}
\labeldef{F} {5\labelsep \labelspace 8}  {nD}
\labeldef{F} {5\labelsep \labelspace 9}  {hD}
\labeldef{F} {5\labelsep \labelspace 10} {ZQ}
\labeldef{L} {5\labelsep \labelspace 9}  {BR}
\labeldef{F} {5\labelsep \labelspace 11} {BCR}
\labeldef{F} {5\labelsep \labelspace 12} {Itwo}
\labeldef{F} {5\labelsep \labelspace 13} {B2}
\labeldef{F} {5\labelsep \labelspace 14} {Vmn}
\labeldef{F} {5\labelsep \labelspace 15} {In}
\labeldef{F} {5\labelsep \labelspace 16} {Rn}

\document

\PaperA4

\hfuzz 10pt

\center
\hrule height 0pt
\vskm.1:.7>

{\twelvepoint\bf \bls1.2\bls
Duality for \KZv/ and \ifMag\\\fi Dynamical Equations
\par}
\vskm1.5:1.4>
\=
\VT/$^{\,\star}$ \ and \ \Varch/$^{\,*}$
\vskm1.5:1.4>
{\it
$^\star$\home/\\ \homeaddr/
\vsk.3>
$^*$\UNC/\\ \UNCaddr/
\vskm1.6:1.5>
\sl June \,2001}
\endcenter

\ftext{\=\bls11pt
$\]^\star\)$\support/ RFFI grant 99\)\~\)01\~\)00101
\>and \,INTAS grant 99\)\~\)01705\vv.06>\nl
\hp{$^*$}{\tenpoint\sl E-mail\/{\rm:} \homemail/}\vv.1>\nl
${\]^*\)}$\support/ NSF grant DMS\)\~\)9801582\vv.06>\nl
\vv-1.2>
\hp{$^*$}{\tenpoint\sl E-mail\/{\rm:} \avemail/}}

\vskm1.7:1.5>

{\ifMag\ninepoint\fi
\Abstract
We consider the \KZv/ (\KZ/\)) and dynamical \eq/s, both \difl/ and \dif/, in
the context of the $\glkn$ duality. We show that the \KZ/ and dynamical \eq/s
naturally exchange under the duality.
\endAbs}

\Sno -1

\vskm1.5:1.1>
\vsk0>

\Sno -1
\sect{Introduction}
The \KZv/ (\KZ/\)) \eq/s is a holonomic system of \difl/ \eq/s for correlation
\fn/s in conformal field theory on the sphere. They play an important role in
\rep/ theory of affine Lie algebras and \qg/s. The \dif/ analogue of the \KZ/
\eq/s is called the quantized \KZv/ (\qKZ/\)) \eq/s.
\par
There are \rat/, \tri/ and elliptic versions of \KZ/ and \qKZ/ \eq/s, depending
on what kind of coefficient \fn/s the \eq/s have. In this paper we consider
the \rat/ and \tri/ versions.
\vsk.1>
In \cite{FMTV} we discovered another holonomic system of \difl/ \eq/s with
\rat/ coefficients, which is compatible with an extension of the \rat/ \KZ/
\eq/s by an element of the Cartan subalgebra. We call the new \eq/s the \rat/
\em{dynamical \difl/} (\DD/\)) \eq/s. The \dif/ analogue of the \DD/ \eq/s was
constructed in \cite{TV1} and called the \rat/ \em{\dif/ dynamical} (\qDD/\))
\eq/s. The \rat/ \qDD/ \eq/s are compatible with the \tri/ \KZ/ \eq/s.
\par
Around 1995 for a simple Lie algebra $\gg$ De~Concini and Procesi introduced
a \$U(\gg)$-valued connection on the set of regular elements of the Cartan
subalgebra of $\gg$ \cite{CP}\). They conjectured that the monodromy of
the connection is described in terms of the corresponding quantum Weyl group.
In the recent work of Toledano\>-Laredo this conjecture is proved for
$\gg=\sln$ \cite{TL}\).
\vsk.1>
When the work \cite{FMTV} had been written, P\]\&Etingof told us about the
De~Concini\>--\)Procesi connection and indicated that their connection is
a special case of our \DD/ \eq/s.
\vsk.1>
The $\glkn$ duality plays an important role in the representation theory
and the classical invariant theory, see \cite{Zh1}\), \cite{H}\).
In \cite{TL} Toledano\>-Laredo discovered that \wrt/ this duality the
De~Concini\>--\)Procesi connection is dual to the \rat/ \difl/ \KZ/ \eq/s
and used this fact to compute the monodromy of the De~Concini\>--\)Procesi
connection in terms of the quantum Weyl group. Our paper is inspired by
the result of \cite{TL}\).
\par
In this paper we show that under the $\glkn$ duality the \KZ/ and \qKZ/ \eq/s
for the Lie algebra $\glk$ exchange with the \DD/ and \qDD/ \eq/s for the Lie
algebra $\gln$, while the \DD/ and \qDD/ \eq/s for $\glk$ exchange with
the \KZ/ and \qKZ/ \eq/s for $\gln$.
\par
The duality between the \KZ/ and \DD/ \eq/s in the \rat/ \difl/ case is
essentially the ``quantum'' version of the duality for isomonodromic
deformation systems \cite{Ha1}\). We thank J\&Harnad for pointing out this
observation to us. The relation of the \difl/ \KZ/ \eq/s and the isomonodromic
deformation systems is described in \cite{R}\), \cite{Ha2}\).
\vsk.1>
An interplay between extremal projectors and \Rms/, relevant to the duality
studied in our paper, has been observed in \cite{ST}\).
\vsk.1>
We consider three cases in the paper: the \rat/ \difl/, \tri/ \difl/, and \rat/
\dif/ cases. To consider the \tri/ \dif/ case, the duality for the \tri/ \qKZ/
\eq/s and the \tri/ \qDD/ \eq/s, one has to employ the \$q$-analogue of
the $\glkn$ duality: the $\Uqkn$ duality. This will be done elsewhere.
The $\Uqkn$ duality has been described in \cite{B}\), \cite{TL}\).
\par
The paper is organized as follows. After introducing basic notation we
subsequently describe the \difl/ \KZ/ and \DD/ \eq/s, and the \rat/ \dif/ \qKZ/
and \qDD/ \eq/s, for the Lie algebra $\glk$. This is done in Sections~2\,--\,4.
The main results are Theorems~\[dualKZDD] and \[BR] in Section~5.
\vsk.2>
The first author thanks Maxim Nazarov for helpful discussions.

\Sect{Basic notation}
Let $n$ be a nonnegative integer. \em{A \part/ $\la=(\la_1\),\la_2\),\ldots{})$
with at most $\>k$ parts} is an infinite nonincreasing sequence of nonnegative
integers \st/ $\la_{k+1}=0$. Denote by $\Pck$ the set of \part/s with at most
$k$ parts and by $\Pc$ the set of all \part/s. We often make use of the
embedding $\Pck\]\to\Ck$ given by truncating the zero tail of a partition:
$(\la_1\)\lc\la_k\),0\>,0\>,\ldots{})\map(\la_1\)\lc\la_k)$. Since obviously
$\Pcm\]\sub\Pck$ for $m\le k$, in fact, we have a collection of embeddings
$\Pcm\to\Ck$ for any $m\le k$. What particular embedding is used will be clear
from the context.
\vsk.1>
Let $e_{ab}$, $a,\bk$, be the standard basis of the Lie algebra $\glk$:
$[\)e_{ab}\>,e_{cd}\)]\)=\)\dl_{bc}\>e_{ad}-\dl_{ad}\>e_{cb}$. We take
the Cartan subalgebra $\hg\sub\glk$ spanned by $e_{11}\lc e_{kk}$, and
the nilpotent subalgebras $\ngp$ and $\ngm$ spanned by the elements $e_{ab}$
for $a<b$ and $a>b$, \resp/. We have the standard Gauss decomposition
$\glk=\ngp\]\oplus\)\hg\oplus\ngm$.
\vsk.1>
Let $\eps_1\lc\eps_k$ be the basis of $\)\hg^*\}$ dual to $e_{11}\lc e_{kk}$:
$\bra\)\eps_a\>,e_{bb}\)\ket\)=\)\dl_{ab}$. We identify $\hg^*\}$ with $\Ck\}$
mapping $\la_1\)\eps_1\lsym+\la_k\>\eps_k$ to $(\la_1\lc\la_k)$.
The root vectors of $\glk$ are $e_{ab}$ for $a\ne b$, the corresponding root
being equal to $\al_{ab}=\)\eps_a\]-\eps_b$. The roots $\al_{ab}$ for $a<b$
are positive. The simple roots are $\alk$: $\al_a=\)\eps_a\]-\eps_{a+1}$.
Denote by $\rho$ a half-sum of positive roots.
\goodbm
\vsk.1>
We choose the standard invariant bilinear form $(\,{,}\,)$ on $\glk$:
$(\)e_{ab}\>,e_{cd}\))\)=\)\dl_{ad}\>\dl_{bc}$. It defines an \iso/
$\hg\to\hg^*\}$. The induced bilinear form on $\hg^*$ is
$(\)\eps_a\>,\eps_b\))\)=\)\dl_{ab}$.
\vsk.1>
The element $h\)=\)e_{11}\lsym+e_{kk}\in\hg$ is central in $\glk$. We identify
\vv.05>
the Lie algebra $\slk$ with the Lie subalgebra of $\glk$ orthogonal to $h$.
Then ${\glk=\slk\oplus\C\)h}$.
\vsk.1>
For a \gkmod/ $W\}$ and a weight $\la\in\hg^*\}$
let $\Wla$ be the weight subspace of $W\}$ of weight $\la$.
\vsk.2>
For any $\la\in\Pck$ we denote by $V_\la$ the \irr/ \gkmod/ with \hw/ $\la$.
By abuse of notation, for any $l\in\Zp$ we write $V_l$ instead of
$V_{(\)l,\)0\)\lc 0)}$. Thus, $V_0=\C$ is the trivial \gkmod/, $V_1=\Ck$ with
the natural action of $\glk$, and $V_l$ is the \$l$-th symmetric power of
$V_1$.
\vsk.1>
The element \,$I\)=\!\sum_{a,b\)=1}^k e_{ab}\>e_{ba}$ is central in $\Uk$.
In the \irr/ \gkmod/ $V_\la$ it acts as multiplication by $(\la\>,\la+2\rho)$.
\vsk.1>
Define a \$\glk$-action on the \pol/ ring $\C\)[\)x_1\lc x_k]$ by \difl/
operators: $e_{ab}\map x_a\der_b$, where $\der_b=\der/\der x_b$, and denote
the obtained \gkmod/ by $\Vb$. Then
$$
\Vb\,=\,\Plus_{l=0}^\8\,V_l\,,
\Tag{Vb}
$$
the submodule $V_l$ being spanned by homogeneous \pol/s of degree $l$.
The \hwv/ of the submodule $V_l$ is $x_1^{\>l}$.

\Sect{\KZv/ and \difl/ dynamical \{\mline operators}
For any $g\in\Uk$ set
$g\'i\:=\>\id\)\lox{}\%{g}_{\Clap{\sss\text{\=$i$-th}}}{}\]\lox\)\id
\in\bigl(U(\glk)\bigr)\vpb{\ox n}\!$.
\vv-.2>
We consider $\Uk$ as a subalgebra of $\Ukn\}$, the embedding \>$\Uk\)\hto\Ukn$
\vvn.06>
being given by the \$n\)$-fold coproduct, that is, $x\>\map x\'1\lsym+x\'n$
for any $x\in\glk$.
\vsk.1>
Let \,$\Om\,=\}\sum_{a,b=1}^k e_{ab}\)\ox\)e_{ba}$ be the Casimir tensor,
\,and let \,$\dsize\Omo=\,{1\over 2}\>\tsum_{a=1}^k e_{aa}\)\ox\)e_{aa}$,
\vvnn0:.1>
$$
\Omp=\,\Omo+\tsum_{1\le a<b\le k} e_{ab}\)\ox\)e_{ba}\,,\qqq
\Omm=\,\Omo+\tsum_{1\le a<b\le k} e_{ba}\)\ox\)e_{ab}\,,
$$
so that $\Om\>=\>\Omp\!+\>\Omm$.
\vskgood-.5:.7>
Fix a nonzero complex number $\ka$. Consider \difl/ operators $\nablasn\]$ and
\>$\hablasn\]$ with coefficients in $\Ukn\}$ depending on complex \var/s $\zn$
and $\lak$:
\vvnn-.5:-.4>
$$
\gather
\nabla_{\}z_i}(z\);\la)\,=\,\ka\>{\der\over\der z_i\]}
\>-\)\tsum_{a=1}^k \la_a\>({e_{aa}})\'i\:
\>-\)\sum_{\tsize{j=1\atop j\ne i}}^n\>{\Om\'{ij}\over z_i-z_j}\;,
\\
\nn-10>
\ald
\nn16>
\habla_{\}z_i}(z\);\la)\,=\,\ka\)z_i\>{\der\over\der z_i\]}\>-\)
\tsum_{a=1}^k\bigl(\la_a\]-{e_{aa}\over 2}\bigr)\>({e_{aa}})\'i\:
\>-\)\sum_{\tsize{j=1\atop j\ne i}}^n\>
{z_i\>\Omp\'{ij}+\)z_j\>\Omm\'{ij}\over z_i-z_j}\;.
\\
\nnm-1:-6>
\mngood
\endgather
$$
The \difl/ operators $\nablasn$ (resp.~$\hablasn$\}) are called
the \em{\rat/} (resp.~\em{\tri/}) \em{\KZv/} (\KZ/\)) \em{operators}.
The following statements are well known.
\Th{ratKZ}
The operators $\>\nablasn\}$ pairwise commute.
\endpro
\Th{trigKZ}
The operators $\>\hablasn\}$ pairwise commute.
\endpro
The \em{\rat/ \KZ/ \eq/s\,} is a system of \difl/ \eq/s
$$
\nabla_{\}z_i}u\,=\,0\,,\qqq \inn\,,\kern-2em
$$
for a \fn/ $u(\zn\);\lak)$ taking values in an \$n\)$-fold tensor product
of \gkmod/s. The \em{\tri/ \KZ/ \eq/s\,} is a system of \difl/ \eq/s
$$
\habla_{\}z_i}u\,=\,0\,,\qqq \inn\,.\kern-2em
$$
\par
Introduce \difl/ operators $\Dlk\]$ and \>$\Dhk\]$ with coefficients
in $\Ukn\}$ depending on complex \var/s $\zn$ and $\lak$:
\ifMag
\vvn-.2>
$$
\align
D_{\la_a}(z\);\la)\,=\,\ka\>{\der\over\der \la_a\]}
\>-\)\tsum_{i=1}^n z_i\>({e_{aa}})\'i\:
\>-{} & {}\)\sum_{\tsize{b=1\atop b\ne a}}^k\>
{e_{ab}\>e_{ba}-\)e_{aa}\over\la_a\]-\la_b}\;.
\\
\nn3>
\Dh_{\la_a}(z\);\la)\,=\,\ka\>\la_a{\der\over\der \la_a\]}\>+
\>{e_{aa}^{\)2}\over2}\>-\)\tsum_{i=1}^n z_i\>({e_{aa}})\'i\:\>-{} &
\\
\nn5>
{}\)-\)\sum_{b=1}^k\,\sum_{1\le i<j\le n}\!(e_{ab})\'i\)(e_{ba})\'j
\>-{} & {}\)\sum_{\tsize{b=1\atop b\ne a}}^k\)
{\la_b\over\la_a\]-\la_b\]}\,(e_{ab}\>e_{ba}-\)e_{aa})\,.
\\
\cnn-.4>
\endalign
$$
\else
\vvn-.4>
$$
\gather
D_{\la_a}(z\);\la)\,=\,\ka\>{\der\over\der \la_a\]}
\>-\)\tsum_{i=1}^n z_i\>({e_{aa}})\'i\:
\>-\)\sum_{\tsize{b=1\atop b\ne a}}^k\>
{e_{ab}\>e_{ba}-\)e_{aa}\over\la_a\]-\la_b}\;.
\\
\nn4>
\Dh_{\la_a}(z\);\la)\,=
\,\ka\>\la_a{\der\over\der \la_a\]}\>+\>{e_{aa}^{\)2}\over2}
\>-\)\tsum_{i=1}^n z_i\>({e_{aa}})\'i\:
\>-\)\sum_{b=1}^k\,\sum_{1\le i<j\le n}\!(e_{ab})\'i\)(e_{ba})\'j
\>-\)\sum_{\tsize{b=1\atop b\ne a}}^k\)
{\la_b\over\la_a\]-\la_b\]}\,(e_{ab}\>e_{ba}-\)e_{aa})\,\Rlap{.}
\\
\cnn-.6>
\endgather
$$
\fi
Recall that $e_{ab}\)=\sum_{i=1}^n\>(e_{ab})\'i\:$.
We call the \difl/ operators $\Dlk$ (resp.~$\Dhk$) the \em{\rat/}
(resp.~\em{\tri/}) \em{\difl/ dynamical} (\)\DD/\)) \em{operators}.
\Th{KZDD}
The operators $\>\nablasn$, $\Dlk\}$ pairwise commute.
\endpro
\nt
The statement follows from the same result for the $\slk$ case
in \cite{FMTV}\).
\Par
\Th{trigDD}
\back\cite{TV2}
The operators $\>\Dhk\}$ pairwise commute.
\endpro
Later we will formulate analogues of Theorem \[KZDD] for the \tri/ \KZ/
operators and for the \tri/ \DD/ operators. They involve \dif/ dynamical
operators and \dif/ (quantized) \KZv/ operators.
\vsk.2>
The \em{\rat/ \DD/ \eq/s\,} is a system of \difl/ \eq/s
\vvn.1>
$$
D_{\la_a}u\,=\,0\,,\qqq \ak\,,\kern-2em
\vv.1>
$$
for a \fn/ $u(\zn\);\lak)$ taking values in an \$n\)$-fold tensor product
of \gkmod/s. The \em{\tri/ \DD/ \eq/s\,} is a system of \difl/ \eq/s
\vvn.1>
$$
\Dh_{\la_a}u\,=\,0\,,\qqq \ak\,.\kern-2em
\vvgood
\mmgood
$$

\Sect{Rational \dif/ dynamical operators}
For any $a\>,\bk$, $a\ne b$, \>introduce a series $B_{ab}(t)$ depending
on a complex \var/ $t$:
\vvn.3>
$$
B_{ab}(t)\,=\,1+\>\sum_{s=1}^\8\>e_{ba}^s\)e_{ab}^s\,
\prod_{j=1}^s\,{1\over j\>(t-\)e_{aa}\]+\)e_{bb}-j)}\;.
\vv.2>
$$
The series has a well-defined action in any \fd/ \gkmod/ $W\}$,
giving an \$\End(W)\)$-valued \raf/ of $t$. These series have zero weight:
$$
\bigl[\)B_{ab}(t)\>,x\)\bigr]\>=\>0\qquad\text{for any}\quad x\in\hg\,,
\vvmm0:-.4>
$$
satisfy the inversion relation
\vvnn-.4:0>
$$
B_{ab}(t)\>B_{ba}(-\)t)\,=\,1\,-{e_{aa}\]-e_{bb}\over t}\;,
\vv-.8>
$$
and the braid relation
\vvnn-.1:-.4>
$$
B_{ab}(t)\>B_{ac}(t+s)\>B_{bc}(s)\,=\,B_{bc}(s)\>B_{ac}(t+s)\>B_{ab}(t)\,.
$$
\Rem
In notation of \cite{TV1} the series $B_{ab}(t)$ equals
$p\)(t-1\);\)e_{aa}\]-e_{bb}\),e_{ab}\),e_{ba})$. The properties of $B_{ab}(t)$
follow from the properties of \fn/s $B_w(\la)$ in the $\slk$ case,
see \cite{TV1\), Section 2.6\)}\).
\enddemo
\Rem
The series $B_{ab}(t)$ first appeared in the definition of the extremal
cocycle in \cite{Zh2}\).
\enddemo
Consider the products $X_1\lc X_k$ depending on complex \var/s $\zn$ and
$\lak$:
\ifMag
\vvn-.7>
$$
\align
X_a(z\);\la)\,= {}& \,\bigl(\)B_{ak}(\la_{ak})\ldots
B_{a,\)a+1}(\la_{a,\)a+1})\)\bigr)\vpb{-1}\x{}
\Tag{Xa}
\\
\nn4>
& {}\x\,\prod_{i=1}^n\,\bigl(z_i^{-\)e_{aa}}\bigr)\'i\,
B_{1a}(\la_{1a}\]-\ka)\ldots B_{a-1,\)a}(\la_{a-1,\)a}\]-\ka)\,,\kern-2em
\endalign
$$
\else
$$
X_a(z\);\la)\,=\,\bigl(\)
B_{ak}(\la_{ak})\ldots B_{a,\)a+1}(\la_{a,\)a+1})\)\bigr)\vpb{-1}
\,\prod_{i=1}^n\,\bigl(z_i^{-\)e_{aa}}\bigr)\'i\,
B_{1a}(\la_{1a}\]-\ka)\ldots B_{a-1,\)a}(\la_{a-1,\)a}\]-\ka)\,,\kern-.3em
\Tag{Xa}
$$
\fi
acting in any \$n\)$-fold tensor product $W_1\lox W_n$ of \fd/ \gkmod/s.
Here $\la_{bc}=\la_b-\la_c$.
\vsk.3>
Denote by $T_u$ a \dif/ operator acting on a \fn/ $f(u)$ by
\ifMag
$$
(T_uf)(u)\,=\,f(u+\ka)\,.
$$
\else
$(T_uf)(u)=f(u+\ka)$.
\fi
Introduce \dif/ operators $\Qlk$:
\vvnn0:-.4>
$$
Q_{\la_a}(z\);\la)\,=\,X_a(z\);\la)\,T_{\la_a}\,.
\vv.2>
$$
We call these operators the (\em{\rat/}\)) \em{\dif/ dynamical} (\qDD/\))
operators.
\Th{KZqDD}
The operators $\>\hablasn$, $\Qlk\}$ pairwise commute.
\endpro
\nt
The statement follows from the same result for the $\slk$ case in \cite{TV1}\).
\goodbm
\Par
\ifMag\else\nt\fi
In more conventional form the equalities
\,$[\)\habla_{\}z_i}\),\)Q_{\la_a}\)]\)=\)0$ \,and
\,$[\)Q_{\la_a}\),\)Q_{\la_b}\)]\)=\)0$ \resp/ look like:
\vvnn-.2:0>
$$
\gather
\habla_{\}z_i}(z\);\la)\,X_a(z\);\la)\,=\,
X_a(z\);\la)\,\habla_{\}z_i}(z\);\la_1\lc\la_a\]+\ka\lc \la_k)\,,
\\
\nnm8:10>
X_a(z\);\la)\,X_b(z\);\la_1\lc\la_a\]+\ka\lc \la_k)\,=\,
X_b(z\);\la)\,X_a(z\);\la_1\lc\la_b\]+\ka\lc \la_k)\,.
\endgather
$$
\par
The \em{\rat/ \qDD/ \eq/s\,} is a system of \deq/s
$$
Q_{\la_a}u\,=\,u\,,\qqq \ak\,,\kern-2em
$$
for a \fn/ $u(\zn\);\lak)$ taking values in an \$n\)$-fold tensor product
of \gkmod/s.

\Sect{Rational \dif/ \KZv/ operators}
Consider the Yangian $\Yk$. The algebra $\Uk$ is embedded in $\Yk$ as
a Hopf subalgebra, and we identify $\Uk$ with the image of this embedding.
\vsk.2>
There is an algebra \hom/ $\ev:\)\Yk\)\to\>\Uk$, \> called the \em{\ehom/},
which is identical on the subalgebra $\Uk\sub\Yk$. The \ehom/ is not a \hom/
of Hopf algebras.
\vsk.2>
The Yangian $\Yk$ has a distinguished one-parametric family of \aut/s $\rho_u$
depending on a complex parameter $u$. For any \gkmod/ $W\}$ we denote by $W(u)$
the pullback of $W\}$ through the \hom/ \>$\ev\o\,\rho_u$. Yangian modules of
this form are called \em{\emod/s}.
\vsk.2>
For any \fd/ \irr/ \gkmod/s $V\),\>W\}$ the tensor products $V(t)\ox W(u)$
and $W(u)\ox V(t)$ are isomorphic \irr/
\vv.1>
\Ykmod/s, provided $t-u\nin\Z$.
The intertwiner can be taken of the form $\PVW\>\RVW(t-u)$ where
\vvmm.1:.08>
\>$\PVW\}:\)V\ox W\to\>W\ox V$ is the flip map:
$\PVW\}:\)v\ox w\)\map\alb\)w\ox v$,
\,and $\RVW(t)$ is a \rat/ \$\End(V\ox W)\)$-valued \fn/,
\vv.1>
called the \em{\rat/ \Rm/} for the tensor product $V\ox W$.
\vsk.2>
We describe below the \Rm/ $\RVW(t)$ in terms of the $\glk$ actions
in $V\}$ and $W\}$. Let $v$ and $w$ be the respective \hwv/s.
The \Rm/ $\RVW(t)$ is uniquely determined by the \$\glk$-invariance,
\vvn.3>
$$
\bigr[\)\RVW(t)\>,\>g\otimes 1+1\otimes g\>\bigr]\,=\,0\qqq
\text{for any}\quad g\in\glk\,,\kern-1.4em
\vv.3>
\Tag{Rinv}
$$
the commutation relations
\vvn-.3>
$$
\RVW(t)\>\bigl(\)t\>e_{ab}\ox 1\,+\tsum_{c=1}^k\>e_{ac}\ox e_{cb}\)\bigr)\,=\,
\bigl(\)t\>e_{ab}\ox 1\,+\tsum_{c=1}^k\>e_{cb}\ox e_{ac}\)\bigr)\>\RVW(t)\,,
\Tag{Rdef}
$$
and the normalization condition
$$
\RVW(t)\,v\ox w\,=\,v\ox w\,.
\vv.2>
$$
The introduced \Rms/ obey the inversion relation
\vvn-.1>
$$
\RVW(t)\,R_{WV}\"{21}(-\)t)\,=\,1
\vv-.2>
\Tag{inv}
$$
where $R_{WV}\"{21}\)=\)P_{WV}\:\>R_{WV}\:\>\PVW$, \,and the \YB/
\vvn.2>
$$
\RUV(t)\>\RUW(t+u)\>\RVW(u)\,=\,\RVW(u)\>\RUW(t+u)\>\RUV(t)\,.
\vv.2>
\Tag{YB}
$$
\par
The aforementioned facts on the Yangian $\Yk$ are well known;
for example, see \cite{MNO}\).
\vsk.2>
Consider the \gkmod/ $\Vb$, and let $v_l=x_1^{\>l}$ be the \hwv/
\vv.06>
of the \irr/ component $V_l\sub\Vb$. We define the \Rm/ $\RVb(t)$
as a direct sum of the \Rms/ $\RVlm(t)$:
$$
\RVb(t)\,v\ox v'\)=\,\RVlm(t)\,v\ox v'\qqq
\text{for any}\quad v\ox v'\]\in\)V_l\ox V_m\,,\kern-2em
$$
the \Rms/ $\RVlm(t)$ being normalized by $\RVlm(t)\,v_l\ox v_m=v_l\ox v_m$.
It is clear that $\RVb(t)$ obeys both the inversion relation and the \YB/.
\vsk.2>
Consider the products $K_1\lc K_n$ depending on complex \var/s $\zn$ and
$\lak$:
\ifMag
\vvn-.6>
$$
\align
\kern.8em K_i(z\);\la)\,= {}& \,\bigl(\)R_{in}(z_{in})\ldots
R_{i,\)i+1}(z_{i,\)i+1})\)\bigr)\vpb{-1}\x{}
\Tag{Ki}
\\
\nn4>
&{}\x\,\prod_{a=1}^k\,\bigl(\)\Rlap{\la_a}{\la}\vpb{-\)e_{aa}}\bigr)\'i\,
R_{1i}(z_{1i}\]-\ka)\ldots R_{i-1,\)i}(z_{i-1,\)i}\]-\ka)\,,\kern-1.4em
\\
\cnn-.1>
\endalign
$$
\else
$$
K_i(z\);\la)\,=\,\bigl(\)R_{in}(z_{in})\ldots
R_{i,\)i+1}(z_{i,\)i+1})\)\bigr)\vpb{-1}
\,\prod_{a=1}^k\,\bigl(\)\Rlap{\la_a}{\la}\vpb{-\)e_{aa}}\bigr)\'i\,
R_{1i}(z_{1i}\]-\ka)\ldots R_{i-1,\)i}(z_{i-1,\)i}\]-\ka)\,,\kern-.3em
\vv-.1>
\Tag{Ki}
$$
\fi
acting in an \$n\)$-fold tensor product ${W_1\lox W_n}$ of \gkmod/s. Here
\ifMag
$z_{ij}=z_i-z_j$, \,and ${R_{ij}(t)\)=\)\bigl(R_{W_iW_j}\:(t)\]\bigr)\'{ij}}$.
\else
${R_{ij}(t)\)=\)\bigl(R_{W_iW_j}\:(t)\]\bigr)\'{ij}}$, \,and $z_{ij}=z_i-z_j$.
\fi
\vsk.2>
Introduce \dif/ operators $\Zzn$:
\vvn-.2>
$$
Z_{z_i}(z\);\la)\,=\,K_i(z\);\la)\,T_{z_i}\,.
$$
We call these operators the (\em{\rat/}\)) \em{quantized} \em{\KZv/} (\qKZ/\))
operators.
\Th{qKZ}
\back\cite{FR}
The operators $\>\Zzn$ pairwise commute.
\endpro
\Th{qKZDD}
\back\cite{TV2}
The operators $\>\Zzn$, $\Dhk\}$ pairwise commute.
\endpro
\nt
The last theorem extends Theorems \[trigDD] and \[qKZ].
\Par
\ifMag\else\nt\fi
In more conventional form the equalities
\,$[\)Z_{z_i}\),\)Z_{z_j}\)]\)=\)0$ \,and
\,$[\)Z_{z_i}\),\)\Dh_{\la_a}\)]\)=\)0$ \resp/ look like:
\vvnn-.3:0>
$$
\gather
K_i(z\);\la)\,K_j(z_1\lc z_i+\ka\lc z_n\);\la)\,=\,
K_j(z\);\la)\,K_i(z_1\lc z_j+\ka\lc z_n\);\la)\,,
\\
\nn8>
\Dh_{\la_a}(z\);\la)\,K_i(z\);\la)\,=\,
K_i(z\);\la)\,\Dh_{\la_a}(z_1\lc z_i+\ka\lc z_n\),\la)\,.
\endgather
$$
\par
The \em{\rat/ \qKZ/ \eq/s\,} is a system of \deq/s
$$
Z_{z_i}u\,=\,u\,,\qqq \inn\,,\kern-2em
\vv-.1>
$$
for a \fn/ $u(\zn\);\lak)$ taking values in an \$n\)$-fold tensor product
of \gkmod/s.
\vskmgood-.8:.8>

\Sect{$\{\glkn$ duality}
In what follows we are going to consider the Lie algebras $\glk$ and $\gln$
simultaneously. In order to distinguish generators, modules, etc., we will
indicated the dependence on $k$ and $n$ explicitly, for example, $e_{ab}\+k\>$,
$V_\la\+n\}$.
\vsk.2>
Consider the \pol/ ring $P_{kn}\)=\)\C\)[\)\xkn\)]$ of $\>k\)n$ \var/s.
We define a \$\glk$-action on $P_{kn}$ by
\vvnn-.3:-.7>
$$
e_{ab}\+k\,\map\)\tsum_{i=1}^n\,x_{ai}\)\der_{b\)i}\,,
\vv-.3>
\Tag{left}
$$
where $\der_{b\)i}=\der/\der x_{b\)i}$,
and we define a \$\gln$-action on $P_{kn}$ by
\vvn-.3>
$$
e_{ij}\+n\,\map\)\tsum_{a=1}^k\,x_{ai}\)\der_{aj}\,.
\vv-.5>
\Tag{right}
$$
There are two natural \iso/s of vector spaces:
$\bigl(\C\)[\)x_1\lc x_k]\)\bigr)\vpb{\ox n}\}\to\,P_{kn}$,
\vvn-.2>
$$
(p_1\lox p_n)\)(\)x_{11}\lc x_{kn})\,=\>
\tprod_{i=1}^n\,p_i(\)x_{1i}\lc x_{ki})\,,
\vv-.4>
\Tag{kiso}
$$
and $\,\bigl(\C\)[\)x_1\lc x_n]\)\bigr)\vpb{\ox k}\}\to\,P_{kn}$,
\vvn-.6>
$$
(p_1\lox p_k)\)(\)x_{11}\lc x_{kn})\,=\>
\tprod_{a=1}^k\,p_a(\)x_{a1}\lc x_{an})\,.
\vv-.2>
\Tag{niso}
$$
We have a simple proposition.
\Prop{Pkn}
The module $P_{kn}$ is isomorphic to $\Vbkn\!$ as a \gkmod/ by \(kiso),
and it is isomorphic to $\Vbnk\!$ as a \gnmod/ by \(niso)\).
\endpro
It is easy to see that the actions \(left) and \(right) commute with
each other, thus making $P_{kn}$ into a module over the direct sum
$\glk\]\oplus\)\gln$. The following theorem is well known, e.g.~see \cite{H}.
\Th{dual}
The $\glk\]\oplus\)\gln$ module $P_{kn}$ has the decomposition
$$
P_{kn}\,=\,\Plus_{\la\)\in\)\Pc\_{\]\min(k,n)}}\ V_\la\+k\]\ox\)V_\la\+n\,.
\vvmm-.8:->
$$
\endpro
Fix vectors $\la=(\)l_1\lc l_n)\in\Zp^{\)n}$ and
\vv-.1>
$\mu=(\)m_1\lc m_k)\in\Zp^{\)k}$ \st/ $\sum_{i=1}^n\)l_i\)=\sum_{a=1}^k\)m_a$.
Let
\vvnn-:-.4>
$$
\Zlm\,=\,\bigl\lb\>
(d_{ai})_{\tsize{a\)=\)1\lc k\atop\Lph ai\>=\)1\lc n}}\]\in\>\Zp^{\)kn}\;
\Vert\tsum_{a=1}^k\)d_{ai}=l_i\,,\quad\tsum_{i=1}^n\)d_{ai}=m_a\)\bigr\rb\,.
\vv-.4>
$$
Denote by $\Plm\sub P_{kn}$ the span of all monomials
$\prod_{a=1}^k\)\prod_{i=1}^n\)x_{ai}^{\)d_{ai}}$ \st/ $(d_{ai})\in\Zlm$.
\vv.3>
The next statement easily follows from Proposition \[Pkn], and formulae
\(Vb)\), \(left) and \(right)\).
\vskmgood-.5:.5>
\Prop{Plm}
The \iso/s \(kiso) and \(niso) induce the \iso/s of the weight subspaces
$\Vlm$ and $\Vml$, and the space $\Plm$.
\endpro
In what follows we will need another description of the above mentioned \iso/s.
It is given below.
\vsk.2>
For an indeterminate $y$ and a nonnegative integer $s$ we define
the \em{divided powers} $y\#s\}$ by the rule: $y\#0\}=\)1$ \,and
\,$y\#s\}=\)y^s\}/\)s!$ \,for $s>0$.
\vsk.2>
Let $v_{l_i}\+k\],\,v_{m_j}\+n$ be \hwv/s of the respective modules
$V_{l_i}\+k\],\,V_{m_j}\+n$.
\Lm{basisk}
A basis of the \wt/ subspace $\Vlm$ is given by vectors
\vv-.3>
$$
v_\db\+k\,=\,\tprod_{a=2}^k\bigl(e_{a1}\+k\bigr)\#{d_{a1}}v_{l_1}\+k\lox
\tprod_{a=2}^k\bigl(e_{a1}\+k\bigr)\#{d_{a\]n}}v_{l_n}\+k\,,
\qqq \db\)=\)(d_{ai})\in\Zlm\,.\kern-2em
\vv-.6>
$$
\endpro
\Lm{basisn}
A basis of the \wt/ subspace $\Vml$ is given by vectors
\vv-.2>
$$
v_\db\+n\,=\,\tprod_{i=2}^n\bigl(e_{i1}\+n\bigr)\#{d_{1i}}v_{m_1}\+n\lox
\tprod_{i=2}^n\bigl(e_{i1}\+n\bigr)\#{d_{ki}}v_{m_k}\+n\,,
\qqq \db\)=\)(d_{ai})\in\Zlm\,.\kern-2em
\vv-.4>
$$
\endpro
For any \,$\db\in\Zlm$ \,set
\,$x\#{\db\)}=\)\prod_{a=1}^k\)\prod_{i=1}^n\)x_{ai}\#{d_{ai}}\}$.
\,Define \iso/s of vector spaces:
\vvn-.4>
$$
\alignat2
& \Vlm\,\to\,\Plm\,,\kern 4em && v_\db\+k\map\>x\#{\db\)},
\Tag{vkx}
\\
\nn8>
& \Vml\,\to\,\Plm\,,\kern 4em && v_\db\+n\map\>x\#{\db\)},
\Tag{vnx}
\\
\cnn-.1>
\endalignat
$$
and extend them by linearity to \iso/s of vector spaces
\,$\bigl(\)\Vb\+k\)\bigr)\vpb{\ox n}\}\to\)P_{kn}$ \,and
$\bigl(\)\Vb\+n\)\bigr)\vpb{\ox k}\}\to\)P_{kn}$.
\Prop{PV}
The maps \(vkx) and \(vnx) coincide with the maps \(kiso) and \(niso), \resp/.
\endpro
\nt
The proof is straightforward.
\Par
Consider the action of \KZ/, \qKZ/, \DD/ and \qDD/ operators for the Lie
algebras $\glk$ and $\gln$ on \$P_{kn}$-valued \fn/s of $\zn$ and $\lak$,
treating the space $P_{kn}$ as a tensor product $\Vbkn\!$ of \gkmod/s,
\vv.06>
and as a tensor product $\Vbnk\!$ of \gnmod/s. We will write $F\simeq G$ if
the operators $F$ and $G$ act on the \$P_{kn}$-valued \fn/s in the same way.
\Lm{Ikn}
We have \,$I\+k\!-\>k\)\sum_{a=1}^k\)e\+k_{aa}\,\simeq
\,I\+n\!-\>n\)\sum_{i=1}^n\)e\+n_{ii}$.
\endpro
\nt
The proof is straightforward.
\Par
Set
\vvn-.5>
$$
C\+k_{ab}(t)\,=\;{\Gm(t+1)\,\Gm(t-e\+k_{aa}\]+e\+k_{bb})\over
\Gm(t-e\+k_{aa})\,\Gm(t+e\+k_{bb}+1)}\;, \kerm 1.5:4>
C\+n_{ij}(t)\,=\;{\Gm(t+1)\,\Gm(t-e\+n_{ii}\]+e\+n_{jj})\over
\Gm(t-e\+n_{ii})\,\Gm(t+e\+n_{jj}+1)}\;.\kerm-1:0>
\mmgood
\vv.3>
\Tag{C}
$$
\Th{dualKZDD}
For any \,$\inn$ \>and \>$\ak$ \>we have
\vvn-.1>
$$
\alignat2
& \qquad\nabla\+k_{\}z_i}(z\);\la)\,\simeq\,D\+n_{z_i}(\la\);z)\,, &
D\+k_{\la_a}(z\);\la)\,\simeq\,\nabla\+n_{\]\la_a}(\la\);z)\,, \kerm-.3:0> &
\Tag{nD}
\\
\nn8>
& \qquad\habla\+k_{\}z_i}(z\);\la)\,\simeq\,\Dh\+n_{z_i}(\la\);z)\,, &
\Dh\+k_{\la_a}(z\);\la)\,\simeq\,\habla\+n_{\]\la_a}(\la\);z)\,, \kerm-.3:0> &
\Tag{hD}
\\
\nn8>
& Z\+k_{z_i}(z\);\la)\,\simeq\,N\+n_i(z)\>Q\+n_{z_i}(\la\);z)\,, \kerm3:3.3> &
N\+k_a(\la)\>Q\+k_{\la_a}(z\);\la)\,\simeq\,Z\+n_{\la_a}(\la\);z)\,.
\ifMag\kern-2.3em\else\kern-2em\fi &
\Tag{ZQ}
\endalignat
$$
Here
\vv->
$$
\align
N\+n_i(z)\, &{}=\]\prod_{1\le j<i}C\+n_{ji}(z_{ji}-\ka)
\}\prod_{i<j\le n}\bigl(C\+n_{ij}(z_{ij})\]\bigr)\vpb{-1}
\\
\nn-4>
\Text{and}
\nn.4>
N\+k_a(\la)\, &{}=\]\prod_{1\le b<a}C\+k_{ba}(\la_{ba}-\ka)
\}\prod_{a<b\le k}\bigl(C\+k_{ab}(\la_{ab})\]\bigr)\vpb{-1}.
\\
\cnn-.5>
\nngood
\endalign
$$
\endpro
\Pf.
Equalities \(nD) and \(hD) for \difl/ operators are verified in
a straightforward way. Equalities \(ZQ) for \dif/ operators follow from
Theorem \[BR], and formulae \(Xa) and \(Ki)\).
\epf
\Rem
The first of equalities \(nD) at $\la=0$ (\resp/, the second one at $z=0$)
has been discovered in \cite{TL}\).
\enddemo
\Th{BR}
For any \)$a\),\bk$, $a\ne b$, \>and any \>$i\),\jn$, $i\ne j$, we have
\vvnn-.4:0>
$$
B\+k_{ab}(t)\>C\+k_{ab}(t)\,\simeq\,R\+n_{ab}(t)\,,
\ifMag\kern3.3em\else\kern4em\fi
R\+k_{ij}(t)\,\simeq\,B\+n_{ij}(t)\>C\+n_{ij}(t)\,.
\vv-.4>
\Tag{BCR}
$$
\endpro
\Rem
Cf.~the observation made after Theorem~6.6 in \cite{FV}\).
\enddemo
\Pf of Theorem \[BR]\).
The Lie algebras $\glk$ and $\gln$ appear in our consideration on equal
footing. So, it suffices to prove only the first formula in \(BCR)\).
Moreover, since both the action of $B\+k_{ab}(t)\>C\+k_{ab}(t)$ on $P_{kn}$
and that of $R\+n_{ab}(t)$ involves only the \var/s $x_{a1}\lc x_{an}$,
$x_{b1}\lc x_{bn}$, it is enough to prove the claim for $k=2$.
\vsk.2>
Let $\la=(l_1\lc l_n)$. Set $|\)\la\)|=l_1\lsym+l_n$.
Consider the decomposition of the \gtmod/ $\Vlt$ into a direct sum of \irr/s:
$$
\Vlt\,=\Plus_{0\le\)m\)\le|\la|/\Rlap{2}}\VtWm\>.
\vv-.1>
$$
Here \>$W_0\), W_1\),\ldots{}$ are multiplicity spaces. The summands $\VtWm$
are distinguished by \eva/s of the central element $I\+2\}$:
\vvnn.6:0>
$$
I\+2\vst{\)\VtWm}\)=\,|\)\la\)|^2\]+\)|\)\la\)|\)+\)2\)m\)(m-|\)\la\)|-1)\,.
\Tag{Itwo}
$$
Set
\vvnn0:-.2>
$$
\Vlmtm\)=\,\bigl(V\+2\'{|\la|-\)m\),\)m}\ox W\+2_m\bigr)\>\cap\>\Vlmt\,.
\vv.4>
$$
It follows from \cite{TV1\), Section 2.5\)} that
\vvn-.6>
$$
B\+2_{12}(t)\vst{\)\Vlmtm}\!=\,
\prod_{j=m}^{m_2-1}\,{t+m_2-j\over t-m_1+j}\;.
\Tag{B2}
$$
\par
Consider the decomposition of the \gnmod/ $\>\Vmn$ into a direct sum of \irr/s:
\vvnm-.4>
$$
\Vmn\,=\Plus_{m\)=0}^{\min\)(m_1,\)m_2)}\Vmmn\,.
\vv-.2>
\Tag{Vmn}
$$
The summands $\Vmmn$
are distinguished by \eva/s of the central element $I\+n\}$:
\ifMag
$$
\align
& I\+n\vst{\)\Vmmn}\!={}
\Tag{In}
\\
\nn6>
& \}{}=\,(m_1\]+m_2)^2\]+\)(n-1)(m_1\]+m_2)\)+\)
2\)m\)(m-m_1\]-m_2\]-1)\,.
\endalign
$$
\else
$$
I\+n\vst{\)\Vmmn}\!=\,(m_1\]+m_2)^2\]+\)(n-1)\)(m_1\]+m_2)\)+\)
2\)m\)(m-m_1\]-m_2\]-1)\,.
\Tag{In}
$$
\fi
Since the \Rm/ $R\+n(u)$ is \$\gln$-invariant, see \(Rinv)\), it acts as
a scalar on each summand of the decomposition \(Vmn). The \eva/s of $R\+n(t)$
can be found from \(Rdef)\), see \cite{KRS}, \cite{KR}\):
\vvnm-.3>
$$
R\+n(t)\vst{\)\Vmmn}\!=\,\prod_{j=0}^{m-1}\,{t-m_1+j\over t+m_2-j}\;.
\vv-.5>
\Tag{Rn}
$$
Set
\vvn-.1>
$$
\Vmltm\)=\,\Vmmn\>\cap\>\Vmlt\,.
$$
\vsk.3>
Let $\mu=(m_1\>,m_2)$. Lemma \[Ikn] and formulae \(Itwo)\), \(In) imply that
under \iso/s \(vkx) and \(vnx) the images of $\Vlmtm$ and $\Vmltm$ in the space
$\Plm$ coincide. Comparing now formulae \(C)\), \(B2) and \(Rn)\), we get
the required statement.
\vvgood
\epf

\myRefs
\widest{FTV}
\parskip.1\bls

\ref\Key B
\by P\]\&Baumann
\paper \$q$-\}Weyl group of a \$q$-Schur algebra \jour Preprint \yr 1999
%% \pages 1\~\)24
\endref

\ref\Key CP
\by C\&De~Concini and C\&Procesi \yr 1995
\paperinfo unpublished
\endref

\ref\Key EV
\by \Etingof/ and \Varch/
\paper Dynamical Weyl groups and applications
\jour Adv\.Math. \vol 167 \yr 2002 \issue 1 \pages 74\>\~127
\endref

\ref\Key FMTV
\by \Feld/, Ya\&Markov, \VT/ and  \Varch/
\paper Differential \eq/s compatible with \KZ/ \eq/s
\jour Math\. Phys., Analysis and Geometry \vol 3 \yr 2000 \pages 139\>\~177
\endref

\ref\Key FR
\by \Fre/ and \Reshy/
\paper Quantum affine algebras and holonomic \dif/ \eq/s
\jour \CMP/ \vol 146 \yr 1992 \pages 1\~\>60
\endref

\ref\Key FV
\by \Feld/ and  \Varch/
\paper Three formulas for eigenfunctions of integrable Schr\"odinger operators
\jour Compositio Math. \vol 107 \yr 1997 \pages 143\)\~175
\endref

\ref\Key H
\by R\& Howe
\paper Perspectives on invariant theory: Schur duality, multiplicity-free
actions and beyond
\jour Israel Math\. Conf\. Proc. \vol 8 \yr 1995 \pages 1\~182
\endref

\ref\Key Ha1
\by J\&Harnad
\paper Dual isomonodromic deformations and moment maps to loop algebras
\jour \CMP/ \vol 166 \yr 1994 \issue 2 \pages 337\)\~\)366
\endref

\ref\Key Ha2
\by J\&Harnad
Quantum isomonodromic deformations and the \KZv/ \eq/s
\jour CRM Proc. Lecture Notes \vol 9 \yr 1996 \pages 155\>\~161
\endref

\ref\Key KR
\by P\]\&P\]\&Kulish and \Reshy/
\paper On \$GL_3$-invariant \sol/s of the \YB/ and associated quantum systems
\jour \JSM/ \vol 34 \yr 1986 \pages 1948\)\~1971
\endref

\ref\Key KRS
\by P\]\&P\]\&Kulish, \Reshy/ and \Skl/
\paper \YB/ and \rep/ theory I
\jour \LMP/ \vol 5 \yr 1981 \pages 393\)\~\)403
\endref

\ref\Key MNO
\by A\&Molev, \MN/ and G\&Olshanski
\paper Yangians and classical Lie algebras
\jour Russian Math\. Surveys \vol 51 \yr 1996 \pages 205\>\~\)282
\endref

\ref\Key MV
\by Ya\&Markov and  \Varch/
\paper Hypergeometric \sol/s of trigonometric \KZ/ \eq/s satisfy dynamical
\deq/s \jour Adv\.Math. \vol 166 \yr 2002 \issue 1 \pages 100\>\~147
\endref

\ref\Key R
\by \Reshet/
\paper The \KZv/ system as a deformation of the isomonodromy problem
\jour \LMP/ \vol 26 \yr 1992 \issue 3 \pages 167\)\~177
\endref

\ref\Key ST
\by Yu\&F\]\&Smirnov and V\]\&N\&\]Tolstoy
\paper Extremal projectors for usual, super and quantum algebras
and their use for solving Yang-Baxter problem
\inbook in Selected topics in QFT and Mathematical Physics,
(Proc.\ of the V\]-th Intern\.Conf., Liblice, Czechoslovakia, June 25\>\~\)30,
1989) \yr 1990 \publ \WSa/
\endref

\ref\Key TL
\by V\]\&Toledano\>-Laredo
\paper A Kohno-Drinfeld theorem for quantum Weyl groups
\jour \DMJ/ \vol 112 \year 2002 \issue 3 \pages 421\~\)451
\endref

\ref\Key TV1
\by \VT/ and \Varch/
\paper Difference \eq/s compatible with \tri/ \KZ/ \difl/ \eq/s
\jour \IMRN/ \yr 2000 \issue 15 \pages 801\~\)829.
\endref

\ref\Key TV2
\by \VT/ and \Varch/
\paper The dynamical \difl/ \eq/s compatible with the \rat/ \qKZ/ \eq/s
\paperinfo in preparation
\endref

\ref\Key Zh1
\by D\&P\]\&Zhelobenko
\book Compact Lie groups and their representations
\bookinfo Transl\. Math\. Mono., vol\point 40 \yr 1983 \publ \AMSa/
\endref

\ref\Key Zh2
\by \Zhel/
\paper Extremal cocycles on Weyl groups \jour \FAA/ \vol 21 \yr 1987
\pages 183\)\~192
\endref

\ref
\by \Zhel/
\paper Extremal projectors and generalized Mickelsson algebras on reductive
Lie algebras \jour Math\. USSR, Izvestiya \vol 33 \yr 1989 \pages 85\>\~100
\endref

\endRefs

\bye